\documentclass{article}

\begin{document}
\newtheorem{cor}{Corollary}
\newtheorem{theo}[cor]{Theorem}
\newtheorem{prop}[cor]{Proposition}
\newtheorem{lem}[cor]{Lemma}
\newtheorem{re}[cor]{Remark}
\newtheorem{ex}[cor]{Example}

\newenvironment{rem} {\begin{re} \upshape} {\end{re}}
\newenvironment{exa} {\begin{ex} \upshape} {\end{ex}}
\newenvironment{proof}{\noindent \emph{Proof.}}  {\vspace{.2cm}}
\renewcommand{\today}{\normalsize  Department of Mathematics \\
University of California \\ Berkeley, CA 94720-3840 \\
\texttt{king@math.berkeley.edu} }


\title{A mass formula for unimodular lattices with no roots}
\author{Oliver D. King}

\maketitle

\begin{abstract}
We derive a mass formula for $n$-dimensional unimodular lattices
having any prescribed root system. We use Katsurada's formula for
the Fourier coefficients of Siegel Eisenstein series to compute
these masses for all root systems of even unimodular
32-dimensional lattices and odd unimodular lattices of dimension
$n\leq 30$. In particular, we find the mass of even unimodular
32-dimensional lattices with no roots, and the mass of odd
unimodular lattices with no roots in dimension $n\leq 30$,
verifying Bacher and Venkov's enumerations in dimensions 27 and
28. We also compute better lower bounds on the number of
inequivalent unimodular lattices in dimensions 26 to 30 than those
afforded by the Minkowski-Siegel mass constants.
\end{abstract}


\section{Introduction}
First we review some definitions. For more information, see
\cite{CS1} and \cite{Kit1}.

An $n$-dimensional \emph{lattice} $\Lambda$ is the set
$\mbox{Z}v_1 + \cdots + \mbox{Z}v_n$ of all integer linear
combinations of a basis $\{v_1, \ldots, v_n\}$ for $\mbox{R}^n.$
We associate with $\Lambda$ the \emph{Gram matrix} $A$ with
$i,j$-th entry the inner product $(v_i,v_j)$. (The matrix $A$ is a
positive definite quadratic form, and much of what follows may be
reformulated in the language of quadratic forms.) The
\emph{determinant} of $\Lambda$ is defined to be $\det(A)$, and
$\Lambda$ is called \emph{unimodular} if its determinant is 1. We
say that $\Lambda$ is \emph{integral} if $(v,v')$ is an integer
for all $v,v' \in \Lambda$. Such a lattice is called \emph{even}
(or \emph{Type II}) if $(v,v)$ is always even, and \emph{odd} (or
\emph{Type I}) otherwise. The \emph{dual} of $\Lambda$ is
$\Lambda' = \{v \in \mbox{R}^n : (v,x) \in \mbox{Z} \mbox{ for all
} x \in \Lambda \}$, and the determinant of $\Lambda$ is the order
of the finite abelian group $\Lambda'/\Lambda$. Since we shall
(with the exception of the occasional dual lattice) primarily be
concerned with integral lattices, we generally omit the adjective
integral in what follows. The \emph{norm} of a vector $v \in
\Lambda$ is defined to be $(v,v)$, the square of its length. We
say that a lattice is \emph{decomposable} (or \emph{reducible}) if
it can be written as an orthogonal direct sum of two nonzero
sublattices, and is \emph{indecomposable} (or \emph{irreducible})
otherwise.

Each $n$-dimensional unimodular lattice $\Lambda$ has a vector $u$
such that $(u,x) \equiv (x,x) \pmod 2$ for all $x \in \Lambda$.
Such a vector is called a \emph{parity vector} \cite[Preface to
3rd edition, p.\ xxxiv]{CS1}, or a \emph{characteristic vector}
\cite{Bo1}, or a \emph{canonical element} \cite{Se}. The set of
parity vectors forms a coset $u + 2\Lambda$ in $\Lambda /
2\Lambda$, and each parity vector satisfies $(u,u) \equiv n \pmod
8$.

Let $N$ be an $n$-dimensional lattice with Gram matrix $A$ and $M$
be an $m$-dimensional lattice with Gram matrix $B$. We say that
$N$ is \emph{represented} by $M$ if there exists an $m \times n$
integral matrix $X$ for which $X^t B X = A$. We let $r(M,N)$
denote the number of representations of $N$ by $M$. When $m=n$, we
say that $M$ and $N$ are (integrally) \emph{equivalent} if there
exists an integral matrix $X$ with determinant $\pm1$ for which
$X^t B X = A$. $M$ and $N$ are in the same \emph{genus} if they
are equivalent over the $p$-adic integers $\mbox{Z}_p$ for each
prime $p$ (including $p=\infty$, for which $\mbox{Z}_p =
\mbox{Q}$). We define $\mbox{Aut}(N)$ to be the group of $n \times
n$ integral matrices $X$ for which $X^t A X = A$. (Note that the
definitions of $\det(N)$, $r(M,N)$ and equivalence are independent
of the integral bases chosen for the lattices, and so is
$\mbox{Aut}(N),$ up to isomorphism.) The \emph{theta series} of a
lattice $\Lambda$ is defined by $\Theta_{\Lambda}(q) = \sum_{v \in
\Lambda} q^{(v,v)} = \sum_{k=0}^{\infty} r(\Lambda,k)q^k$ where
$q=e^{\pi iz}.$

For our purposes, a \emph{root} is a vector of norm 1 or 2, and
the \emph{root system} of $\Lambda$ is the set of roots in
$\Lambda$. The lattice generated by a root system is called a
\emph{root lattice}, and we define the \emph{rank} of a root
system to be the dimension of the corresponding root lattice. Root
lattices are completely classified: they are direct sums of the
irreducible root lattices $\mbox{Z}$, $A_n\, (n\geq1), D_n\, (n
\geq 4), E_6, E_7, \mbox{ and } E_8$. We use the same notation to
refer to root systems, and for brevity we sometimes write the root
system $n_1 A_1 \oplus \cdots \oplus n_j E_8$ as $A_1^{n_1}\cdots
E_8^{n_j}.$

In this paper we are concerned with the problem of classifying
unimodular lattices, and also the subproblem of classifying
unimodular lattices without roots (which correspond to denser
sphere-packings), up to equivalence.

The Minkowski-Siegel mass formula (see \cite{CS3}) gives the sum
of the reciprocals of the orders of the automorphism groups of all
inequivalent lattices in a given genus. The mass constants can be
used to verify that an enumeration of inequivalent lattices in a
given genus is complete. They also give a lower bound for the
number of inequivalent lattices in a given genus (sometimes called
the \emph{class number}).

In dimensions divisible by 8 there are two genera of unimodular
lattices, Type I and Type II (except in dimension 0, in which
there is just one lattice of Type II). In all other dimensions
there is only one genus of unimodular lattices, Type I. The lower
bounds provided by the mass constants show that the number of
unimodular lattices increases super-exponentially as a function of
the dimension. Unimodular lattices have been completely enumerated
in dimensions $n \leq 25$, but number more than 900 million in
dimension 30. Unimodular lattices without roots have been
completely enumerated in dimensions $n \leq 28$, but number more
than $8 \times 10^{20}$ in dimension 33. In both cases, somewhere
between dimension 26 and dimension 32 there is a transition from
being completely classified to being numerous enough to make
classification unappealing, so information about what is going on
in these dimensions is of interest.

Our approach is based on a suggestion by Borcherds that the
Fourier coefficients of Siegel Eisenstein series of degree $4k$
could be used to derive something analogous to the
Minkowski-Siegel mass formula, but which gives the mass of all the
even unimodular lattices of dimension $8k$ having any given root
system.

We wrote a computer program which uses Katsurada's formula (see
\cite{Kat1}) for the Fourier coefficients to calculate the masses
for all possible root systems of even unimodular 32-dimensional
lattices. From these masses, we used the methods of \cite[Chapter
16]{CS1} to compute the masses, for each root system, of
unimodular lattices in dimensions $n \leq 30$.

In particular, we have a mass formula for those lattices which
have the empty root system (that is, which have no roots).  This
formula verifies the known masses in dimensions $n \leq 28$, and
provides new lower bounds for the number of unimodular lattices
without roots in dimensions 29 to 32.

We also used our program to compute better lower bounds on the
total number of odd unimodular lattices in dimensions 26 to 30,
and even unimodular lattices in dimension 32, than those gotten
from the Minkowski-Siegel mass constants.

Our results may be viewed as a coarse classification of lattices,
as inequivalent lattices may have the same root system. For
dimensions $n \leq 23$, in which it happens to be the case that a
unimodular lattice in a given genus is completely determined by
its root system, our results coincide with the previously known
enumerations. The same is true for the even unimodular
24-dimensional lattices (which are known as the \emph{Niemeier
lattices}).

\section{A mass formula for even unimodular lattices having any
given root system} \label{matrixform}

 Let $\Omega$ be the set of inequivalent even unimodular lattices of
 dimension $8k$.
 We define the \emph{mass} $m$ of $\Omega$ by $$m = \sum_{\Lambda \in \Omega}
 \frac{1}{|\mbox{Aut}(\Lambda)|}.$$

By the Minkowski-Siegel mass formula, for $k>0$ we have
 $$m  =
\frac{|B_{4k}|}{8k} \prod_{j=1}^{4k-1} \frac{|B_{2j}|}{4j} $$
where $B_i$ is the $i$th Bernoulli number. (See \cite{CS1} and
\cite{CS3} for this, and for the mass formulae for other genera of
lattices.)

For an $n$-dimensional lattice $N$, we define a weighted average
number of representations of $N$ by the lattices $\Lambda \in
\Omega$ by $$a(N) = \frac{1}{m}\sum_{\Lambda \in \Omega} \frac{r(
\Lambda,N)} {|\mbox{Aut}(\Lambda)|}.$$

Let $\{R_1, ..., R_s\}$ be the set of all the root lattices of
dimension $n$ or less with no vectors of norm 1. Each $R_i$ is the
direct sum of lattices $A_j (j\geq 1), D_j (j \geq 4), E_6,
E_7,\mbox{ and } E_8$. Let $\Omega_i$ be the set of lattices in
$\Omega$ having root system $R_i$, and let $m(R_i) = \sum_{\Lambda
\in \Omega_i} |\mbox{Aut}(\Lambda)|^{-1}$ be the mass of those
lattices in $\Omega$ which have root system $R_i$, so that $m =
m(R_1) + \cdots + m(R_s)$.

\begin{prop} Let $U$ be the $s \times s$ matrix with
$i,j$-th entry $r(R_j, R_i)$, let $v$ be the vector $(m(R_1),
\ldots, m(R_s))^t$, and let $w$ be the vector $(a(R_1), \ldots,
a(R_s))^t$. Then  $\frac{1}{m}U v = w$. Furthermore, $U$ is
invertible, so $m(R_i) = (m U^{-1}w)_i$ gives the mass of the
lattices in $\Omega$ having root system $R_i$.
\end{prop}

\begin{proof} Observe that if $R$ is any root
lattice and if $S$ is the root system of a lattice $\Lambda$ then
$r(\Lambda,R) = r(S,R)$, so for $j=1, \ldots, s$ we have
\begin{eqnarray*}
a(R_j) & = & \frac{1}{m}\sum_{\Lambda \in \Omega}
\frac{r(\Lambda,R_j)}{|\mbox{Aut}(\Lambda)|} = \frac{1}{m}
\sum_{i=1}^{s} \sum_{\Lambda \in \Omega_i}
\frac{r(\Lambda,R_j)}{|\mbox{Aut}(\Lambda)|} \\ & = & \frac{1}{m}
\sum_{i=1}^{s} r(R_i,R_j) \sum_{\Lambda \in \Omega_i}
\frac{1}{|\mbox{Aut}(\Lambda)|} \\ & = & \frac{1}{m}
\sum_{i=1}^{s}r(R_i,R_j)\,m(R_i).
\end{eqnarray*}

Thus  $\frac{1}{m}U v = w$.  We may assume the $R_i$'s are ordered
so that their dimensions are non-decreasing, and so that among
those with the same dimension the determinants are non-increasing.
With this ordering, $r(R_j, R_i) = 0$ whenever $i>j$, so $U$ is
upper triangular. Since each diagonal element $r(R_i, R_i) =
|\mbox{Aut}(R_i)|$ is positive, we have $\det(U) \neq 0$ so $v = m
U^{-1}w$.
\end{proof}

\begin{rem} The values $a(R_j)$ are the Fourier coefficients of Siegel
Eisenstein series. Borcherds, Freitag and Weissauer \cite{BFW}
used a relation similar to $\frac{1}{m}U v = w$ to compute the
coefficients of a cusp form from the known values of $m(R_i)$ in
dimension 24. We shall do the inverse in dimension 32: use the
values of $a(R_j)$ to derive the values $m(R_i)$. We discuss how
to compute $a(R_j)$ in Section~\ref{computeA}. The problem of
computing $r(R_i,R_j)$ is largely combinatorial; we discuss it in
Section~\ref{computeR}.
\end{rem}

\section{Masses of 32-dimensional even unimodular lattices with any
given root system}

Let $m_n^{\mbox{\scriptsize II \normalsize}}(R)$ and
$m_n^{\mbox{\scriptsize I \normalsize}}(R)$ denote the masses of
the $n$-dimensional even and odd unimodular lattices having root
system $R$, and let $m_n(R)$ denote their sum. Let $w(R)$ denote
the order of the Weyl group of $R$; $w(R)$ is the product of the
orders of the Weyl groups of the irreducible components of $R$,
where $w(A_n) = (n+1)!, w(D_n) = 2^{n-1} n! , w(E_6) = 2^7 3^4 5
\cdot 7, w(E_7) = 2^{10} 3^4 5 \cdot 7$, and  $w(E_8) = 2^{14} 3^5
5^2 7$. It is sometimes more convenient to list the values of
$m_n^{\mbox{\scriptsize II \normalsize}}(R) \cdot w(R)$ than it is
to list the values of $m_n^{\mbox{\scriptsize II \normalsize}}(R)$
alone, but the latter can easily be recovered from the former.

We used a computer to calculate $m_{32}^{\mbox{\scriptsize II
\normalsize}}(R)$ for each root system $R$ with $\mbox{rank}(R)
\leq 32$. The computation took about two weeks on a Sun Ultra 60,
running a program written in Common Lisp and compiled with Franz
Inc.'s Allegro CL.\@ (We discuss several issues related to the
implementation in Sections~\ref{computeA} to \ref{hacks}.)

\begin{theo} The nonzero values of $m_{32}^{\mbox{\scriptsize II
\normalsize}}(R) \cdot w(R)$ with $\mbox{rank}(R) \leq 9$ are as
listed in Table~\ref{somemasses}. (The list of all $13218$ nonzero
values of $m_{32}^{\mbox{\scriptsize II \normalsize}}(R) \cdot
w(R)$ is available electronically at \cite{Kin}.)
\end{theo}

\renewcommand{\arraystretch}{1.3} \normalsize

\begin{table}
\caption{Masses of 32-dimensional even unimodular lattices}
\label{somemasses} \small
\begin{center}
 \(
\begin{array}{rccrc}
\dim & \mbox{root system } R & m_{32}^{\mbox{\scriptsize II
\normalsize}}(R) \cdot w(R) & \mbox{(as decimal)} &
\mbox{comments}\\  \hline  0 & \emptyset &
\frac{1310037331282023326658917}{238863431761920000 } & 5484461.50
&
  \mbox{no roots} \\ 1
& A_1  & \frac{111536168182433}{5677056}  & 19646832.00 &  \\ 2 &
A_1^2 & \frac{72024731351193941}{1857945600}  & 38765792.00 &
\\ 2 & A_2  & \frac{1327104974887}{2939328}  &    451499.44 &  \\ 3 & A_1^3
& \frac{6904800898075}{124416}  &  55497696.00 &  \\ 3 & A_1 A_2 &
\frac{977951251237}{445440}  &   2195472.50 &  \\ 3 & A_3  &
\frac{329127961}{74240}  &      4433.30 &  \\ 4 & A_1^4  &
\frac{30223371257980501}{471859200}  &  64051672.00 &  \\ 4 &
A_1^2 A_2  & \frac{19867101805}{3456}  &   5748582.50 &  \\ 4 &
A_2^2  &  \frac{1772535692573}{42598400}  &     41610.38 &  \\ 4 &
A_1 A_3 &  \frac{21073837}{768}  &     27439.89 &  \\ 4 & A_4  &
 \frac{8397751}{384000}  & 21.87 &  \\ 4 & D_4  &
 \frac{35841940559}{157212057600}
&         0.23 & \mbox{see \cite{BV}}
\\ 5 & A_1^5  &  \frac{14457125482723}{230400}  &  62747940.00 &  \\ 5 &
A_1^3 A_2  &  \frac{10626384230783}{995328}  &  10676264.00 &  \\
5 & A_1 A_2^2 &  \frac{673556587}{2560}  &    263108.06 &  \\ 5 &
A_1^2 A_3 &  \frac{200386803709}{2211840}  &     90597.34 &  \\ 5
& A_2 A_3 &  \frac{8085187}{5760}  &      1403.68 &  \\ 5 & A_1
A_4 &
 \frac{46917823}{269568} & 174.05 &  \\ 5 & A_1 D_4  &  \frac{473}{240}  & 1.97 &
\\ 5 & A_5 &  \frac{73}{960}  &         0.08 &  \\ 5 & D_5  &
\frac{433}{3317760} & 0.00 & \mbox{see \cite{BV}}
\\ 6 & A_1^6  &  \frac{355695555290333}{6635520}  &  53604776.00 &
\\ 6 & A_1^4 A_2  &  \frac{22484458507}{1440}  &  15614208.00 &  \\ 6 & A_1^2
A_2^2 &  \frac{16820220686833}{19169280}  &    877457.06 &  \\ 6 &
A_1^3 A_3 &  \frac{5452147363}{25920}  &    210345.19 &  \\ 6 &
A_2^3 &  \frac{9785018477}{1866240}  &      5243.17 &  \\ 6 & A_1
A_2 A_3 &  \frac{680633479}{61440}  &     11078.02 &  \\ 6 & A_1^2
A_4 &  \frac{11264777}{15360}  &       733.38 &  \\ 6 & A_3^2  &
 \frac{323400013}{20447232}  &        15.82 &  \\ 6 & A_1^2 D_4  &
 \frac{622763}{69120}  &         9.01 &  \\ 6 & A_2 A_4  &  \frac{6059}{512}  &
11.83 &  \\ 6 & A_2 D_4  &  \frac{344077}{2446080}  &         0.14
&
\\ 6 & A_1 A_5 &  \frac{123853}{159744}  &         0.78 &  \\ 6 & A_1 D_5
&  \frac{1}{648} & 0.00 &
\\ 6 & A_6  &  \frac{1}{4608}  &         0.00 &  \\
\\
\end{array}
\)
\end{center}
\end{table}

\addtocounter{table}{-1}

\begin{table}

\caption{Masses of 32-dimensional even unimodular lattices
(cont.)} \label{moremasses} \small
\begin{center}
 \(
\begin{array}{rccrc}
\dim & \mbox{root system } R & m_{32}^{\mbox{\scriptsize II
\normalsize}}(R) \cdot w(R) & \mbox{(as decimal)} &
\mbox{comments} \\

 \hline

  6 & D_6  &  \frac{1}{18720000}  &
0.00 &  \mbox{unique} \\
  6 & E_6  &  \frac{1}{1268047872}  & 0.00 &
\mbox{unique} \\
 7 & A_1^7 &  \frac{23307759701}{576}  &  40464860.00 &  \\
 7 & A_1^5 A_2
&  \frac{290429642677}{15360}  &  18908180.00 &  \\ 7 & A_1^3
A_2^2 &  \frac{32712111919}{16128}  &   2028280.70 &  \\ 7 & A_1^4
A_3  &
 \frac{491983248817}{1290240}  &    381311.40 &  \\ 7 & A_1 A_2^3  &
 \frac{641243179}{15360}  &     41747.60 &  \\ 7 & A_1^2 A_2 A_3  &
 \frac{528426689}{11520}  &     45870.37 &  \\ 7 & A_1^3 A_4  &
 \frac{828763}{384}
& 2158.24 &  \\ 7 & A_2^2 A_3  &  \frac{238819303}{552960}  &
431.89 &
\\
7 & A_1^3 D_4  &  \frac{2242333}{77760}  &        28.84 &  \\ 7 &
A_1 A_3^2 &  \frac{2750969}{17920}  &       153.51 &  \\ 7 & A_1
A_2 A_4 &  \frac{57178267}{483840}  &       118.18 &  \\ 7 & A_1^2
A_5 & \frac{24071}{5760} & 4.18 &  \\ 7 & A_1 A_2 D_4  &
\frac{959}{640} & 1.50 &
\\ 7 & A_3 A_4 &  \frac{297043}{860160}  &         0.35 &  \\ 7 & A_2 A_5
&  \frac{51383}{774144} & 0.07 &  \\ 7 & A_3 D_4  &
\frac{133}{30720} & 0.00 &
\\ 7 & A_1^2 D_5  &  \frac{11}{1152}  &         0.01 &  \\ 7 & A_1 A_6  &
 \frac{521}{188160}  & 0.00 &  \\ 7 & A_2 D_5  & \frac{1}{6912}  &         0.00 &
\\ 7 & A_7  &  \frac{1}{1376256}  &         0.00 &  \\ 8 & A_1^8  &
 \frac{2005621383142854931}{73987522560}  &  27107562.00 &  \\ 8 & A_1^6
A_2 &                            \frac{2015000372681}{103680}  &
19434804.00 &
\\ 8 & A_1^4 A_2^2  &
\frac{1327413084613}{368640} &   3600838.50 &  \\ 8 & A_1^5 A_3  &
\frac{980176289}{1728}  & 567231.70 &  \\ 8 & A_1^2 A_2^3  &
\frac{13761649541}{80640}  & 170655.37 &  \\ 8 & A_1^3 A_2 A_3  &
\frac{42164571593}{322560}  & 130718.54 &  \\ 8 & A_2^4 &
\frac{176324027322323}{180592312320} & 976.37 &  \\ 8 & A_1^4 A_4
& \frac{9542951321}{1935360}  & 4930.84 &
\\ 8 & A_1 A_2^2 A_3 &  \frac{691208023}{161280}  &      4285.76 &  \\ 8 &
A_1^4 D_4  &  \frac{371231029}{5160960}  &        71.93 &  \\ 8 &
A_1^2 A_3^2  &  \frac{6071138573}{7741440}  &       784.24 &  \\ 8
& A_1^2 A_2 A_4  &  \frac{1655917}{2688}  &       616.04 &  \\ 8 &
A_2 A_3^2 &  \frac{1794539}{120960} & 14.84 &  \\ 8 & A_1^3 A_5  &
 \frac{2245489}{143360}  & 15.66 &
\\

\end{array}
\)
\end{center}
\end{table}

\addtocounter{table}{-1}

\begin{table}

\caption{Masses of 32-dimensional even unimodular lattices (cont.)
} \label{evenmoremasses} \small
\begin{center}
 \(
\begin{array}{rccrc} \dim & \mbox{root system } R & m_{32}^{\mbox{\scriptsize II
\normalsize}}(R) \cdot w(R) & \mbox{(as decimal)} &
\mbox{comments} \\

 \hline
8 & A_1^2 A_2 D_4  &  \frac{6113}{720}  &         8.49 &  \\ 8 &
A_2^2 A_4 &
 \frac{237}{40}  &         5.93 &  \\
 8 & A_1 A_3 A_4  &
\frac{200995}{48384} & 4.15 &  \\ 8 & A_2^2 D_4  &
\frac{18917}{266112}  & 0.07 &  \\ 8 & A_1 A_2 A_5  &
\frac{7739}{8960}  &         0.86 &  \\
  8 & A_1^3 D_5  &
\frac{829}{20160}  & 0.04 &  \\ 8 & A_1 A_3 D_4  &
\frac{353}{6720} & 0.05 &
\\ 8 & A_1^2 A_6  & \frac{2287}{120960}  &         0.02 &  \\ 8 & A_4^2  &
\frac{432673039}{170311680000}  &         0.00 &  \\ 8 & A_3 A_5 &
\frac{1049}{483840}  &         0.00 &  \\ 8 & A_1 A_2 D_5  &
\frac{43}{20160}  & 0.00 &  \\ 8 & A_2 A_6  & \frac{1}{3780}  &
0.00 &
\\ 8 & A_4 D_4 & \frac{1}{15360}  &         0.00 &  \\ 8 & D_4^2 &
\frac{1867}{1937768448}  & 0.00 &  \\ 8 & A_3 D_5  &
\frac{19}{1935360}  & 0.00 &
\\ 8 & A_1^2 D_6  & \frac{607}{85155840}  &         0.00 &  \\ 8 & A_1 A_7
& \frac{5}{1064448}  & 0.00 &  \\ 8 & A_8  & \frac{1}{185794560}
& 0.00 &
\\ 8 & A_2 E_6 & \frac{1}{277136640}  &         0.00 & \mbox{unique} \\ 8 & D_8  &
\frac{1}{1002795171840}  & 0.00 & \mbox{odd Leech} \\ 8 & E_8  &
\frac{1}{8315553613086720000} & 0.00 &  \mbox{Leech}\\

 9 & A_1^9  &                           \frac{4668288705497}{290304}  &  16080690.00 &  \\
9 & A_1^7 A_2  & \frac{5526347655971}{322560}  & 17132774.00 &  \\
9 & A_1^5 A_2^2  & \frac{740728763}{144}  &   5143950.00 &
\\ 9 & A_1^6 A_3  & \frac{1568559349553}{2211840}  &    709164.94 &  \\ 9 &
A_1^3 A_2^3  & \frac{1368567316381}{2903040}  &    471425.60 &  \\
9 & A_1^4 A_2 A_3  & \frac{1271443287}{4480}  &    283804.30 &  \\
9 & A_1 A_2^4 & \frac{23738023}{2560}  &      9272.67 &  \\ 9 &
A_1^5 A_4  & \frac{2142275851}{233280}  &      9183.28 &  \\ 9 &
A_1^2 A_2^2 A_3 & \frac{927969995}{43008}  &     21576.68 &  \\ 9
& A_1^5 D_4  & \frac{158461}{1080} & 146.72 &  \\ 9 & A_1^3 A_3^2
& \frac{36859871}{13440}  & 2742.55 &
\\ 9 & A_1^3 A_2 A_4  & \frac{237492443}{107520}  &      2208.82 &  \\ 9 &
A_2^3 A_3 & \frac{3078967}{17280}  &       178.18 &  \\ 9 & A_1^4
A_5  & \frac{457601}{10080} & 45.40 &  \\ 9 & A_1^3 A_2 D_4  &
\frac{21454849}{645120}  & 33.26 &  \\

\end{array}
\)
\end{center}
\end{table}

\addtocounter{table}{-1}

\begin{table}

\caption{Masses of 32-dimensional even unimodular lattices
(cont.)} \label{evenevenmoremasses} \small
\begin{center}
 \(
\begin{array}{rccrc} \dim & \mbox{root system } R & m_{32}^{\mbox{\scriptsize II
\normalsize}}(R) \cdot w(R) & \mbox{(as decimal)} &
\mbox{comments} \\

 \hline

 9 & A_1 A_2 A_3^2  &
\frac{6031423}{32256} &       186.99 &
\\ 9 & A_1 A_2^2 A_4  & \frac{989231}{13440}  &        73.60 &  \\

 9 & A_1^2
A_3 A_4 & \frac{1223851}{46080} & 26.56 &  \\ 9 & A_1^2 A_2 A_5  &
\frac{6696269}{1182720} & 5.66 &  \\ 9 & A_1 A_2^2 D_4  &
\frac{2927}{2880} & 1.02 &
\\

 9 & A_3^3 & \frac{31211779}{163296000}  &         0.19 &  \\ 9 & A_1^4 D_5
& \frac{1029287}{7464960}  &         0.14 &  \\ 9 & A_1^2 A_3 D_4
& \frac{66941}{184320}  &         0.36 &  \\ 9 & A_2 A_3 A_4  &
\frac{3827}{8064}  & 0.47 &  \\ 9 & A_1^3 A_6  &
\frac{34339}{384000} & 0.09 &
\\ 9 & A_2^2 A_5 & \frac{783383}{13063680}  &         0.06 &  \\ 9 & A_1
A_4^2  & \frac{1021}{32256} & 0.03 &  \\ 9 & A_2 A_3 D_4  &
\frac{1}{252} & 0.00 &
\\ 9 & A_1^2 A_2 D_5 & \frac{37}{2520}  &         0.01 &  \\ 9 & A_1 A_3 A_5  &
\frac{4813}{129024} & 0.04 &  \\ 9 & A_1 A_2 A_6  &
\frac{547}{126720} & 0.00 &
\\ 9 & A_1 A_4 D_4 & \frac{1}{1260}  & 0.00 &  \\ 9 & A_2^2 D_5  & \frac{1}{8000}  &
0.00 &
\\ 9 & A_1^2 A_7 & \frac{43}{430080}  &         0.00 &  \\ 9 & A_1 A_3 D_5 &
\frac{1}{15840} & 0.00 &  \\ 9 & A_4 A_5  & \frac{13}{295680}  &
0.00 &
\\ 9 & A_3 A_6 & \frac{1}{466560}  &         0.00 &  \\ 9 & A_5 D_4  &
\frac{1}{1512000} & 0.00 & \\ 9 & A_2 A_7  & \frac{1}{504000}  &
0.00 &
\\ 9 & A_1^3 E_6 & \frac{1}{1512000} &         0.00 & \mbox{unique} \\ 9 & A_1 A_2 D_6  &
\frac{1}{645120}  & 0.00 & \\ 9 & A_4 D_5 & \frac{1}{1774080}  &
0.00 &
\\ 9 & A_1 A_8  & \frac{1}{20401920} & 0.00 & \\ 9 & D_4 D_5  &
\frac{1}{41287680} & 0.00 & \\ 9 & A_3 D_6  & \frac{1}{52254720} &
0.00 &  \\ 9 & A_1^2 D_7 & \frac{1}{177408000}  & 0.00 &
\mbox{unique}
\\ 9 & A_3 E_6  & \frac{1}{489646080}  & 0.00 & \mbox{unique} \\ 9 & A_9  & \frac{1}{3592512000}  &
0.00 & \\ 9 & A_2 E_7  & \frac{1}{991533312000}  & 0.00 &
\mbox{unique}
\\ 9 & D_9  & \frac{1}{84610842624000}  & 0.00 & \mbox{shorter
Leech} \\
\end{array}
\)
\end{center}
\end{table}

\renewcommand{\arraystretch}{1.0} \normalsize
\normalsize


The root system $R$ of a lattice $\Lambda$ is called
\emph{complete} if the sublattice of $\Lambda$ generated by $R$
has finite index in $\Lambda$, or equivalently if $\mbox{rank}(R)
= \dim(\Lambda)$. The classification of even unimodular lattices
with complete root systems is closely related to the
classification of certain self-dual codes (see \cite{Ve}).

\begin{cor}[Kervaire \cite{Ke}] There are $119$ complete root systems $R$
which occur as root systems of indecomposable $32$-dimensional
even unimodular lattices.
\end{cor}

\begin{proof} There are 143 root systems $R$ for
which $\mbox{rank}(R) = 32$ and $m_{32}^{\mbox{\scriptsize II
\normalsize}}(R) > 0$. Of these, the root system $D_{16}^2$ and
the 23 root systems containing $E_8$ correspond to decomposable
lattices. The remaining 119 root systems correspond to
indecomposable lattices.
\end{proof}

Kervaire \cite{Ke}, extending the work of Koch and Venkov
\cite{KV}, has proven a stronger result: there are exactly 132
indecomposable even unimodular 32-dimensional lattices with
complete root systems, with 119 different root systems occurring.

\begin{cor} There are at least $10000000$ $32$-dimensional even
unimodular lattices without roots.\end{cor}

\begin{proof}
Each such lattice has at least two automorphisms, so the number of
such lattices is at least $2 \times m_{32}^{\mbox{\scriptsize II
\normalsize}}(\emptyset) = 1.096 \times 10^{7}.$ We will have more
to say about lattices without roots is Section~\ref{secnoroots}.
\end{proof}

\begin{rem}
If there is a single even unimodular 32-dimensional lattice with
root system $R$ then the corresponding value
$m_{32}^{\mbox{\scriptsize II \normalsize}}(R) \cdot w(R)$ is of
the form $1/q$ for some $q \in \mbox{N}$; the converse often holds
(but not always: for $R = A_1^7 A_2^3 A_3 D_7$,
$m_{32}^{\mbox{\scriptsize II \normalsize}}(R) \cdot w(R)$ splits
as $1/4 = 1/12 + 1/6$; see \cite{Bo1}). We give one example below,
and shall see more examples in Section~\ref{sec30} (Niemeier
lattices) and Section~\ref{secnoroots} (unimodular lattices with
no roots). The last column of Table~\ref{somemasses} labels as
``unique'' those lattices whose uniqueness follows from
\cite{Bo1}, and in some cases gives the name of the lower
dimensional lattice which when glued to $R$ produces the unique
even unimodular 32-dimensional lattice with root system $R$.
\end{rem}

\begin{exa} In \cite{EG}, Elkies and Gross construct a 26-dimensional
even lattice $L_0$ of determinant 3 with no roots, for which
$|\mbox{Aut}(L_0)| = 2^{13} 3^5 7^2 13 = 1268047872$, and refer to
a preprint \cite{El1} containing a proof of its uniqueness which
uses Euclidean lattices. (Borcherds  \cite[Chapter 5.7]{Bo1} had
previously proved its existence and uniqueness using Lorentzian
lattices). By \cite{Bo1}, even 26-dimensional lattices of
determinant 3 with no roots are in one-to-one correspondence with
even unimodular 32-dimensional lattices having root system $E_6$,
where the order of the automorphism group of the 32-dimensional
lattice is $w(E_6)$ times the order the automorphism group of the
corresponding 26-dimensional lattice. Since
$m_{32}^{\mbox{\scriptsize II \normalsize}}(E_6) \cdot w(E_6) =
|\mbox{Aut}(L_0)|^{-1}$, $L_0$ must be the unique even
26-dimensional lattice of determinant 3 with no roots, and there
must also be a unique 32-dimensional even unimodular lattice with
root system $E_6$. (This also follows from the classification of
27-dimensional unimodular lattices with no roots in \cite{BV}, as
the one-to-one correspondence mentioned above also extends to
27-dimensional unimodular lattices with no roots and with a parity
vector of norm 3 \cite{Bo1}.)
\end{exa}

\begin{rem} In an earlier draft we pointed out the large mass of
lattices with root systems $A_1^k$ for small $k$, and remarked
that roots seem to have a propensity for being orthogonal to one
another. Peters has since sent us a preprint \cite{Pe} explaining
this: he observes that the sum of the masses
$m_{32}^{\mbox{\scriptsize II \normalsize}}(\emptyset) +
\sum_{k=1}^{32} m_{32}^{\mbox{\scriptsize II \normalsize}}(A_1^k)$
in \cite{Kin} is 97.25\% of the total mass of the genus of even
unimodular 32-dimensional lattices, and shows that a lower bound
of 97.11\% of the total mass can be derived from the small size of
the Fourier coefficient corresponding to $A_2$ in the Siegel
Eisenstein series of degree 2 and weight 16.  The idea is that an
even unimodular lattice represents $A_2$ if and only if it has
non-orthogonal roots $u$ and $v$ (with $u \neq \pm v$).

\end{rem}

\section{Masses of unimodular lattices of dimension $n\leq 30$ with
any given root system}

\label{sec30}

In \cite[Chapter 16]{CS1}, Conway and Sloane describe a
correspondence between unimodular lattices of dimension $n\leq 23$
and orbits of norm 4 vectors in even unimodular 24-dimensional
lattices.  They use this correspondence to produce, from the list
of Niemeier lattices, a list of the root systems and automorphism
group orders of all unimodular lattices of dimension $n \leq 23$.
Unimodular lattices of dimension $n\leq 31$ likewise correspond to
orbits of norm 4 vectors in even unimodular 32-dimensional
lattices,  and we can use this correspondence to compute $m_n(R)$
for all root systems $R$ and all $n \leq 30$.

Let $v = 2e$ be a vector of norm 4 in an even unimodular lattice
$\Lambda$ of dimension $n=32$. Then $L_{31} = \{x \in
\Lambda^{\perp}: x+ ne \in \Lambda \mbox{ for some } n \in
\mbox{Z}\}$ is an odd 31-dimensional unimodular lattice, and
$|\mbox{Aut}(L_{31})| = 2|\mbox{Aut}(\Lambda)|/c(v)$, where $c(v)$
is the number of images of $v$ under $\mbox{Aut}(\Lambda)$. If the
lattice $L_{31}$ has exactly $2k$ vectors of norm 1, we can write
$L_{31} = L_{31-r} \oplus \mbox{Z}^{k}$, where $L_{31-k}$ has
minimal norm 2 and $|\mbox{Aut}(L_{31-k})| =
|\mbox{Aut}(L_{31})|/(2^{k} k!)$.  (We shall call $L_{31-k}$ the
\emph{reduced} lattice corresponding to $\Lambda$ and $v$.) This
construction gives a one-to-one correspondence between orbits of
norm 4 vectors in 32-dimensional even unimodular lattices and
unimodular lattices with no vectors of norm 1 in dimensions less
than 32. (See \cite{CS1} or \cite{Bo1}).

We shall mainly be concerned with norm 4 vectors $v$ which are the
sum of two orthogonal roots $r$ and $s$. For such vectors, knowing
just the root system of $\Lambda$ allows us to compute the
dimension and the root system of the reduced lattice $L_{31-k}$.
Some useful information, distilled from \cite[Table 16.8]{CS1}, is
provided in the table below. $$
\begin{array}{lllc@{\hspace{4pt}}ll@{\hspace{10pt}}l}

R & \#r & \hat{R} & \mbox{\ \ shape of } v \mbox{ \ \ } & \#v &
\tilde{R}& \dim(L_{31-k}) \\
\hline
\vspace{-8pt} \\
 A_n & 4  {n+1 \choose 2} & A_{n-2} & \mbox{any}& 6
{n+1 \choose 4} & A_{n-4} & 29 \vspace{1pt}
\\

D_4 & 24 & A_1^3  & \mbox{any}&  24 & \emptyset & 28\\

D_n \,(n \geq 5) & 4  {n \choose 2} & D_{n-2} &
 \multicolumn{4}{l}{\begin{array}{ll@{\hspace{15pt}}ll}(\pm 1^4, 0^{n-4}) & 16
                 {n \choose 4}& D_{n-4}  & 28 \\
                 (\pm 2, 0^{n-1}) & 2n & \emptyset & 32-n
                 \end{array} }  \\

E_6 & 72 & A_5 & \mbox{any}&  270 & \emptyset & 27\\

E_7 & 126 & D_6 & \mbox{any}&  756 & A_1 & 26\\

E_8 & 240 & E_7 & \mbox{any}&   2160 & \emptyset & 24\\

\end{array} $$

The interpretation of this table is as follows:

Case 1: $r$ and $s$ are from different components $R$ and $S$ of
the root system of $\Lambda$. Then the reduced lattice $L_{31-k}$
has dimension 30, and has the same root system as $\Lambda$,
except for the components $R$ and $S$ from which $r$ and $s$ are
taken, which are transformed to $\hat{R}$ and $\hat{S}$ as given
in the third column. The second column gives the number of roots
$\#r$ in each component $R$, from which we can compute the number
of vectors $v = r+s$ with $r \in R$ and $s \in S$.

Case 2: $r$ and $s$ are from the same component $R$ of the root
system of $\Lambda$. Then the dimension of the reduced lattice is
given in last column, and the root system of $L_{32-k}$ is the
same as that of $\Lambda$ except the component $R$ will be
replaced by $\tilde{R}$. Note that for $R = D_n$ with $n \geq 5$,
$\tilde{R}$ and $\dim (L_{32-k})$ depend on the shape of $v$. The
column headed $\#v$ gives the number of norm 4 vectors of that
shape in the component $R$.

\begin{rem} Norm 4 vectors that are not the sum of two orthogonal roots
correspond to reduced lattices of dimension 31. For this reason we
have not computed the mass of unimodular 31-dimensional lattices
with root system $R$. (That mass could be computed from the mass
of even unimodular 40-dimensional lattices with root system $R
\oplus D_9$.)
\end{rem}

Computing the number $c(v)$ (and hence the order of
$\mbox{Aut}(L_{31-k})$ in terms of the order of
$\mbox{Aut}(\Lambda)$)  sometimes requires additional information
about $\Lambda$, such as whether there are automorphisms of
$\Lambda$ which permute the components of its root system with
multiplicity greater than one. But for our purposes this is not an
impediment, as we demonstrate in the following example.

\begin{exa} Let $\Lambda$ be an even unimodular 32-dimensional lattice
with root system $A_1^4 D_5$ and let $v = s + t$ where $s$ is a
root from one of the four components $A_1$ and $t$ is a root from
the component $D_5$. Then the corresponding reduced lattice
$L_{31-k}$ is 30-dimensional, has root system $A_1^4 A_3$ (since
$D_3$ = $A_3$), and has automorphism group order
$|\mbox{Aut}(L_{30})| = (2|\mbox{Aut}(\Lambda)|/c(v))/(2^1 \cdot
1!) = |\mbox{Aut}(\Lambda)|/c(v)$. Since each component $A_1$ has
2 roots and the component $D_5$ has 40 roots, there are $4\cdot 2
\cdot 40 =320$ such vectors $v$, which form anywhere from one to
four orbits under $\mbox{Aut}(\Lambda)$, depending on whether
there are any automorphisms permuting the components $A_1$.
Because of this ambiguity, we do not know exactly what $c(v)$ is.
But suppose the 320 vectors break into $m$ orbits, $V_1, \ldots,
V_m$, with representatives $v_1, \ldots, v_m$. Then
$$\sum_{i=1}^{m} c(v_i) = \sum_{i=1}^{m} |V_i| = 320 $$ so the
total mass of the lattices which correspond to the lattice
$\Lambda$ and any of these vectors $v$ is $$\sum_{i=1}^{m}
\frac{c(v_i)}{|\mbox{Aut}(\Lambda)|} =
\frac{320}{|\mbox{Aut}(\Lambda)|}. $$ Since from
Table~\ref{somemasses} the mass of all 32-dimensional lattices
$\Lambda_j$ with root system $A_1^4 D_5$ is $\sum_{j}
|\mbox{Aut}(\Lambda_{j})|^{-1} = 1029287/7464960$, these lattices
contribute $\sum_{j} 320 |\mbox {Aut}(\Lambda_{j})|^{-1} =
320\cdot 1029287/7464960$ towards the total mass of 30-dimensional
unimodular lattices with root system $A_1^4 A_3$. By similarly
accounting for contributions from orbits of norm 4 vectors in even
unimodular 32-dimensional lattices with other root systems, we can
compute the exact mass of the 30-dimensional unimodular lattices
with root system $A_1^4 A_3$.
\end{exa}

We have, in this manner, computed the mass of $n$-dimensional
unimodular lattices having any given root system for all $n\leq
30.$  We have not provided a table of these masses (as they can
easily be derived from the table of masses of 32-dimensional even
unimodular lattices with any given root system in \cite{Kin}), but
we will use these masses in Section~\ref{classsec} to find lower
bounds on the class numbers of unimodular lattices in dimensions
up to 30.

\begin{exa} Even unimodular 24-dimensional lattices with root
system $R$ correspond to even unimodular 32-dimensional lattices
with root system $R \oplus E_8$, and $m_{24}^{\mbox{\scriptsize II
\normalsize}}(R) = m_{32}^{\mbox{\scriptsize II \normalsize}}(R
\oplus E_8) \cdot w(E_8) \cdot c $ where $c$ is the multiplicity
of $E_8$ in $R \oplus E_8$. We see from the table in \cite{Kin}
that there are 24 root systems of the form $R\oplus E_8$ for which
$m_{32}^{\mbox{\scriptsize II \normalsize}}(R \oplus E_8) > 0$.
Since even unimodular 24-dimensional lattices happen to be
uniquely determined by their root systems, the corresponding
values of $m_{24}^{\mbox{\scriptsize II \normalsize}}(R)$ are
precisely of the form $|\mbox{Aut}(\Lambda)|^{-1}$, where
$\Lambda$ is the Niemeier lattice having root system $R$.
\end{exa}

\section{Mass formula for unimodular lattices with no roots}
\label{secnoroots}

Recall that $m_n(\emptyset)$ denotes the mass of the
$n$-dimensional unimodular lattices with no roots. By considering
all the orbits of norm 4 vectors $v=r+s$ in even unimodular
32-dimensional lattices that correspond to reduced lattices with
no roots, we have \vspace{0pt}
 \begin{eqnarray*}
  m_n(\emptyset) & = & m_{32}^{\mbox{\scriptsize II \normalsize}}(D_{32-n})\, w(D_{32-n})
  \mbox{ for } n \leq 26, n \neq 24\\
  m_{24}(\emptyset) & = & m_{32}^{\mbox{\scriptsize II \normalsize}}(D_8)\, w(D_8)
            + m_{32}^{\mbox{\scriptsize II \normalsize}}(E_8)\, w(E_8) \\
  m_{27}(\emptyset) & = & m_{32}^{\mbox{\scriptsize II \normalsize}}(D_5)\, w(D_5)
            + m_{32}^{\mbox{\scriptsize II \normalsize}}(E_6) \,w(E_6)\\
  m_{28}(\emptyset) & = & m_{32}^{\mbox{\scriptsize II \normalsize}}(D_4)\, 3 w(D_4)
            + m_{32}^{\mbox{\scriptsize II \normalsize}}(D_5)\,w(D_5)\\
  m_{29}(\emptyset) & = & m_{32}^{\mbox{\scriptsize II \normalsize}}(A_3)\, w(A_3)
            + m_{32}^{\mbox{\scriptsize II \normalsize}}(A_4)\,w(A_4)\\
  m_{30}(\emptyset) & = & m_{32}^{\mbox{\scriptsize II \normalsize}}(A_1^2 )\, w(A_1^2)
            + m_{32}^{\mbox{\scriptsize II \normalsize}}(A_1 A_2)\, w(A_1 A_2)
            + m_{32}^{\mbox{\scriptsize II \normalsize}}(A_2^2 )\, w(A_2^2)
  \end{eqnarray*}

These masses $m_n(\emptyset)$ may be computed by looking up the
values of $m_{32}^{\mbox{\scriptsize II \normalsize}}(R)$ in
Table~\ref{somemasses} and \cite{Kin}; the results are listed in
Table~\ref{massnoroots}, split into $m_{n}^{\mbox{\scriptsize I
\normalsize}}(\emptyset)$ and $m_{n}^{\mbox{\scriptsize II
\normalsize}}(\emptyset)$ for $n$ divisible by 8.
Table~\ref{massnoroots} also includes lower bounds on
$m_{n}^{\mbox{\scriptsize I \normalsize}}(\emptyset)$ for $n = 31$
and 32, which we explain below.

\begin{table}
\caption[Masses of $n$-dimensional unimodular lattices without
roots]{Masses of $n$-dimensional unimodular lattices without
roots. Odd lattices appear above, even lattices below.}
\label{massnoroots}
\renewcommand{\arraystretch}{1.4}\normalsize

\begin{center}
\begin{tabular}{cclll}
  $n$ & $m_n^{\mbox{\scriptsize I \normalsize}}(\emptyset)$
   & (as decimal)\,\,\, & \multicolumn{2}{l}{number of lattices} \\
  \hline

  0-22 & 0 & 0 & 0 \\

  $23$ & $\frac{1}{84610842624000}$ & $1.18 \times 10^{-14}$ & 1 &(shorter Leech lattice)\\

  $24$ & $\frac{1}{1002795171840}$ & $9.97 \times 10^{-13}$ & 1  & (odd Leech lattice)\\

  $25$ & 0 & 0 & 0 \\

  $26$ &  $\frac{1}{18720000}$ & $5.34 \times 10^{-8}$ & 1 & (classified in \cite{Bo1})\\

  $27$ & $\frac{206867}{1585059840}$ & $1.31 \times 10^{-4}$ & 3 & (classified in \cite{BV})\\

  $28$ & $\frac{17924389897}{26202009600}$ & $6.84 \times 10^{-1}$ & 38 & (classified in \cite{BV})\\

  $29$ & $\frac{49612728929}{11136000}$ & $4.46 \times 10^3$ &  \multicolumn{2}{l}{more than $8900$} \\

  $30$ & $\frac{7180069576834562839}{175111372800}$ & $ 4.10 \times 10^7$
  & \multicolumn{2}{l}{more than $82000000$} \\

$31$ & ? & $
> 4 \times 10^{11}$ & \multicolumn{2}{l}{ more than $8 \times
10^{11}$ } \\

$32$ & ? & $> 5 \times 10^{15}$ & \multicolumn{2}{l}{more than $1
\times 10^{16} $}\\

$33$ &  ? & $> 4 \times 10^{20}$ &\multicolumn{2}{l}{ more than $8
\times 10^{20}$ \cite{CS2}}

\\
  \\

  $n$ & $m_n^{\mbox{\scriptsize II \normalsize}}(\emptyset)$ & (as decimal)
  &  \multicolumn{2}{l}{number of lattices} \\

   \hline

  $0$ &  1 & 1 & 1  & (empty lattice)\\

  $8$ & 0 & 0 & 0 & \\

  $16$ & 0 & 0 & 0 & \\

  $24$ &  $\frac{1}{8315553613086720000}$ & $ 1.20 \times 10^{-19}$
  & 1 & (Leech lattice)\\

  $32$  & $\frac{1310037331282023326658917}{238863431761920000}$
  & $5.48 \times 10^{6}$ &  \multicolumn{2}{l}{more than $10000000$} \\

\end{tabular}
\end{center}
\renewcommand{\arraystretch}{1.0}\normalsize

\end{table}

Suppose $\Lambda$ is an even unimodular 32-dimensional lattice
with no roots.  It follows from a theta function argument as in
\cite[Theorem 7.17]{CS1} that $\Lambda$ has 146880 vectors $v$ of
norm 4. For each orbit of these vectors, the corresponding reduced
lattice is 31-dimensional and has no roots. Thus these lattices
$\Lambda$ contribute $$\frac{146880}{2} \,
m_{32}^{\mbox{\scriptsize II \normalsize}} (\emptyset)= \frac
{22270634631794396553201589}{55292461056000} \approx 4.03 \times
10^{11}$$ towards the mass of 31-dimensional lattices with no
roots. There are almost certainly additional contributions from
norm 4 glue vectors $v$ in even unimodular 32-dimensional lattices
with root systems $A_1^k$ for some $k$, but we will not attempt to
account for these. We can restate this in terms of parity vectors:


\begin{prop}
The mass of $31$-dimensional unimodular lattices with no roots and
with no parity vectors of norm $7$ is $(146880/2) \,
m_{32}^{\mbox{\scriptsize II \normalsize}} (\emptyset)$.
\end{prop}

\begin{proof}
Let $L$ be a 31-dimensional unimodular lattice with no roots, and
$\Lambda$ the corresponding 32-dimensional even unimodular
lattice. As in \cite[p.\ 414]{CS1}, $\Lambda$ is equal to $\{x+y
\, | \,x \in L^i, y \in \frac{i}{2} + 2\mbox{Z}, \mbox{for }
i=0,1,2,3\} \subset (L^0)' \oplus \frac{1}{2}\mbox{Z}$, where
$L^0$ is the sublattice of $L$ consisting of vectors of even norm;
$L^0, L^1, L^2$ and $L^3$ are the cosets of $L^0$ in its dual
$(L^0)'$; and $L^0 \cup L^2 = L$. The only ways the vector $x+y$
in $\Lambda$ can have norm 2 is if $y=0$ and $(x,x) = 2$ with $x
\in L^0$, if $y=1$ and $(x,x)=1$ with $x \in L^2$, or if $y=\pm
1/2$ and $(x,x)=7/4$ with $x \in L^1$ or $x \in L^3$. In the first
two cases $x$ would be in $L$ which would contradict $L$ having no
roots. Thus $\Lambda$ has a root if and only if there is a vector
$x$ of norm 7/4 in $L^1 \cup L^3$. But in this case $x$ is in
$(L^0)' \setminus L$ --- the so-called \emph{shadow} of $L$ ---
and $2x$ has norm 7. Since the parity vectors of $L$ are precisely
twice the shadow vectors \cite[Preface to 3rd edition, p.\
xxxiv]{CS1}, $L$ has a parity vector of norm 7 if and only if
$\Lambda$ has a root. (The 31-dimensional unimodular lattices with
no roots and no parity vectors of norm 7 must have parity vectors
of norm 15 by \cite{El2}.)
\end{proof}

We can construct 32-dimensional odd unimodular lattices with no
roots from 32-dimensional even unimodular lattices with no roots
as follows: Let $\Lambda$ be a 32-dimensional even unimodular
lattice with no roots, and let $b \in \Lambda / 2 \Lambda$ be a
nonzero element with norm divisible by 4. Then there is a unique
32-dimensional odd lattice $L$ containing $\Lambda_b := \{v \in
\Lambda \, | \, (v,b) \in 2\mbox{Z} \}$ (see \cite[Chapter
0.2]{Bo1}). ($\Lambda$ and $L$ are \emph{neighbors}, meaning their
intersection has index two in each of them.) If $b$ is not
represented by a vector of norm 0 or 4 then $L$ has no roots. Of
the $2^{32}$ elements in $\Lambda / 2 \Lambda$, $2^{31} + 2^{15}$
have norm congruent to $0 \bmod 4$, by Milgram's formula
\cite[Appendix 4]{MH} applied to the lattice $\sqrt{2}\,\Lambda$.
(Milgram's formula says that if $M$ is an even lattice and $M'$
its dual, then $$\sum_{u \in M' / M} \exp(2 \pi i (u,u)/2) = (\det
M)^{1/2} \exp(2 \pi i \sigma /8)$$ where $\sigma$ is the signature
of $M$, which is equal to its dimension when $M$ is positive
definite.) Note that $(b,b) \bmod 4$ depends only on $b \bmod 2
\Lambda$. Vectors $v$ and $w$ of norm 4 in $\Lambda$ are
equivalent mod  $2 \Lambda$ if and only if $v = \pm w$, so
146880/2 elements of $\Lambda / 2 \Lambda$ are represented by
vectors of norm 4. Thus there are $2^{31} + 2^{15} - 146880/2 - 1
= 2147442975$ elements of $\Lambda / 2 \Lambda$ represented by
vectors of norm divisible by 4 but not represented by vectors of
norm 0 or 4. Since each odd 32-dimensional unimodular lattice has
two even neighbors, this implies that
 $m_{32}^{\mbox{\scriptsize I \normalsize}}(\emptyset)
\geq  (2147442975/2) \cdot \, m_{32}^{\mbox{\scriptsize II
\normalsize}} (\emptyset) = 5.89 \times 10^{15}$. (This is almost certainly
an underestimate, since an even neighbor of an odd lattice with no
roots can have roots.)

As an immediate consequence of the values we computed for
$m_n(\emptyset)$ in Table~\ref{massnoroots}, we have:

\begin{cor}[{\cite{CS1,Bo1,BV}}] There exist odd unimodular $n$-dimensional
lattices without roots for $n=23, 24, \mbox { and } 26 - 32$, but
not for $n = 1 - 22$ or $n=25$.
\end{cor}




Unimodular lattices with no roots are known to exist in all
dimensions $n\geq 26$ (see \cite{CS2}, \cite{NS}). In dimensions
$n \leq 28$ they have already been completely enumerated, and so
the masses $m_n(\emptyset)$ may also be computed by summing the
reciprocals of the automorphism group orders of these lattices.
(Happily, this agrees with our mass formula in each case.) For the
even lattices, there is the empty lattice in dimension 0, and the
Leech lattice $\Lambda_{24}$ in dimension 24. The Leech lattice
$\Lambda_{24}$ was discovered by Leech in 1965 \cite{Le} and was
shown to be the unique even unimodular 24-dimensional lattice
without roots by Niemeier \cite{Ni} and by Conway \cite[Chapter
12]{CS1} around 1969. For the odd lattices, there is the shorter
Leech lattice in dimension 23, the odd Leech lattice $O_{24}$ in
dimension 24 \cite{OP}, the lattice $S_{26}$ in dimension 26, 3
lattices in dimension 27, and 38 lattices in dimension 28. The
lattice $S_{26}$ was constructed by Conway in the 1970's and was
shown to be the unique unimodular 26-dimensional lattice without
roots by Borcherds in 1984 \cite{Bo1}. Borcherds also found one of
the 27-dimensional lattices. The full enumerations in dimensions
27 and 28 are due to Bacher and Venkov \cite{BV}.

\begin{rem} In the cases of the Leech lattice and $S_{26}$, with our
mass formula the uniqueness follows immediately from the
constructions, simply by verifying that
$|\mbox{Aut}(\Lambda_{24})|^{-1} = m_{24}^{\mbox{\scriptsize II
\normalsize}}(\emptyset)$ and $|\mbox{Aut}(S_{26})|^{-1} =
m_{26}^{\mbox{\scriptsize I \normalsize}}(\emptyset)$.
\end{rem}

Examples of unimodular lattices with no roots have been
constructed for dimensions 29 to 32 (including the 15 exceptional
even unimodular 32-dimensional lattices classified in \cite{KV}),
and a nonconstructive analytic argument shows that they exist in
all dimensions $n \geq 33$ (see the Conway-Thompson Theorem
\cite[p.\ 46]{MH} for $n \geq 37$, and \cite{CS3} for $33 \leq n
\leq 36$.) In fact for $n \geq 33$ this argument gives a lower
bound for the mass of lattices without roots that is close to the
total mass of the genus, so there are a great many lattices
without roots (see Remark~\ref{bannai}). The idea is that the
coefficients $a_1$ and $a_2$ of the average theta series
$$\frac{1}{m} \sum_{\Lambda \in \Omega}
\frac{\Theta_{\Lambda}(q)}{|\mbox{Aut}(\Lambda)|} = 1 + a_1 q^1 +
a_2 q^2 + a_3 q^3 + \cdots$$ give the average number of vectors of
norm 1 and norm 2, taken over all the lattices in the genus
$\Omega$. If $a_1 + a_2$ is less than 2, then there must be some
$\Lambda \in \Omega$ with no vectors of norm 1 or 2 --- this is
the case when $\Omega$ is the genus of n-dimensional unimodular
lattices for $n \geq 33$. For $n=33$, $a_1+ a_2$ is approximately
1.42, and this implies that $m_{33}^{\mbox{\scriptsize I
\normalsize}}(\emptyset) \geq 4.04 \times 10^{20}$, so there are
more than $8 \times 10^{20}$ 33-dimensional unimodular lattices
without roots \cite{CS2}. For $n \leq 32$, the average number of
roots is greater than 2, so this argument does not apply. Notice
that the coefficient $a_i$ is the same as the average number of
representations $a(N)$ we defined in Section~\ref{matrixform}, in
the special case where $N$ is the $1 \times 1$ matrix $(i)$.

\begin{rem} In dimensions $n \leq 33$, odd unimodular lattices with
no roots have minimal norm 3, except in dimension 32, in which
they can have minimal norm 3 or 4. In dimensions 24 and 32, even
unimodular lattices with no roots have minimal norm 4. See
\cite{CS2}.
\end{rem}

\begin{rem}
\label{bannai} Let $m_n'$ be the mass of odd $n$-dimensional
unimodular lattices with only trivial automorphisms. There are no
such lattices for $n \leq 28$ \cite{BV}, but Bacher has found one
for $n=29$ \cite{Bac}. Bannai \cite{Ban} showed that $m_n' / m_n
\rightarrow 1$ as $n \rightarrow \infty$. For $n>1$ any lattice
with roots has nontrivial automorphisms, so $m'_n /m_n \leq
m_n(\emptyset)/m_n.$  Below we list $m_n(\emptyset)/m_n$ for $26
\leq  n \leq 30$, and lower bounds on $m_n(\emptyset)/m_n$ for
$31\leq n \leq 33$.

\end{rem}

\begin{center}

\begin{tabular}{l@{\hspace{30pt}}l}

\begin{tabular}{cr@{\,}ll}
$n$ & \multicolumn{2}{l}{$m_n(\emptyset)/m_n$}& ref. \\

\hline

26 & & 0.000116 & \cite{Bo1}\\

27 & & 0.000856 & \cite{BV} \\

28 & & 0.00658 &  \cite{BV}\\

29 & & 0.0300 & \\

\end{tabular}  &

\begin{tabular}{cr@{\,}ll}
$n$ & \multicolumn{2}{l}{$m_n(\emptyset)/m_n$}& ref. \\

\hline

30 & & 0.0908 & \\

31 &$>$&0.135 & \\

32 &$>$&0.136 & \\

33 &$>$&0.287 & \cite{CS2} \\

\end{tabular}

\end{tabular}


\end{center}

\section
{Lower bounds on class numbers} \label{classsec}
 Let $\Omega$ be the set of inequivalent lattices in a genus of dimension $n>0$,
 ${m}$ the mass of that genus, and $m(R)$ the mass of those
 lattices having root system $R$.
 Each lattice $\Lambda \in \Omega$ has at least two
 automorphisms, $1 : x \mapsto x$ and $-1 : x \mapsto -x$;
 from this we get the well-known lower bound $|\Omega| \geq \lceil 2 m
 \rceil $. For each  root $r$ of $\Lambda$, the reflection
 $$x \mapsto x - 2 \frac {(x,r)} {(r,r)} r$$ is also in $\mbox{Aut}(\Lambda)$.
 Define $w'(R)$ to be the order of the subgroup of $\mbox{Aut}(\Lambda)$
 generated by reflections and by $-1$. Then $w'(R)$ = $w(R)$ if the map
 $-1$ is already in the Weyl group of $R$, and $w'(R) = 2w(R)$
otherwise. (The Weyl group of $R$ contains $-1$ if and only if
$\mbox{rank}(R) = \dim(\Lambda)$  and each component of $R$ is
$A_1, E_7, E_8$ or $D_k$ for even $k$.) Then there are at least
$\lceil m(R)\, w'(R) \rceil$ lattices with
 root system $R$, so we get an improved lower bound,
 $$|\Omega| \geq \sum_{R} \lceil m(R)\, w'(R) \rceil .$$

 We can do slightly better still, as follows: Write $m(R)\,w'(R) = q
 + a/b$ with $q,a,b \in \mbox{Z}, a<b, \mbox{ and } \gcd(a,b)=1$, and define a
 modified ceiling function by $\langle q + a/b \rangle = q$ if $a=0$,
 $q+1$ if $a=1$, and $q+2$ if $a>1$. It can easily be shown that there are
  at least $\langle m(R)\, w'(R) \rangle$
 lattices with root system $R$, so that
 \begin{equation} \label{eqlrb}
 |\Omega| \geq \sum_{R} \langle m(R)\, w'(R) \rangle .
 \end{equation}

Evaluating this sum with the value of $m_{32}^{\mbox{\scriptsize
II \normalsize}}(R)$ we computed for each $R$ gives $|\Omega| \geq
1162109024.$

\begin{cor} There are at least $1000000000$ even unimodular
$32$-dimensional lattices.
\end{cor}

For comparison, $\lceil 2 m_{32}^{\mbox{\scriptsize II
\normalsize}}\rceil = 80618466$, so our lower bound is 14.4
times larger than the lower bound obtained by doubling the
Minkowski-Siegel mass constant.

We also computed $m_n(R)$ for each $n \leq 30$ and each root
system $R$, and used equation~\ref{eqlrb} to find lower bounds on
the number of odd unimodular lattices. (In the cases where a
single $m_n(R)$ is expressed as a sum of masses corresponding to
different root systems of 32-dimensional even unimodular lattices,
such as  $m_{28}(\emptyset)  =  m_{32}^{\mbox{\scriptsize II
\normalsize}}(D_4)\, 3 w(D_4) + m_{32}^{\mbox{\scriptsize II
\normalsize}}(D_5)\,w(D_5)$, we bounded the number for each
summand individualy, which gives a better overall bound.)

Unimodular lattices in dimensions $n\leq 25$ have been completely
enumerated: even unimodular lattices of dimension 8 by Mordell, of
dimension 16 by Witt, and of dimension 24 by Niemeier \cite{Ni}
(see also Venkov \cite{Ve}); odd unimodular lattices of dimension
$n\leq 16$ by Kneser \cite{Kn}, of dimension $n\leq 23$ by Conway
and Sloane \cite[Chapter 16]{CS1}, and of dimension 24 and 25 by
Borcherds \cite{Bo1}.

Table~\ref{class30} gives our lower bound $\beta_n$ on the number
of unimodular lattices in dimension $n \leq 30$, and our
computation of the number of distinct root systems $r_n$ that
occur in these lattices (including root systems with components
$\mbox{Z}$, which do not occur for even lattices). The table
includes for comparison the actual number $\alpha_n$ of unimodular
lattices in dimension $n\leq 25$, taken from \cite[Table
2.2]{CS1}, and the Minkowski-Siegel mass constants $m_n$, taken
from \cite[Tables 16.3 and 16.5]{CS1}. (The counts of odd
unimodular lattices include those with vectors of norm 1 to
facilitate comparison with $m_n$; since any integral lattice with
a vector of norm 1 is of the form $\mbox{Z}^k \oplus \Lambda$
where $\Lambda$ has minimal (nonzero) norm 2, the counts of
lattices with no vectors of norm 1 can be recovered from
$\beta_n$; similarly, the number of distinct root systems with no
components $\mbox{Z}$ can be recovered from $r_n$, with the
caution that 8 of the 24 root systems of even unimodular
24-dimensional lattices also occur as root systems of odd
unimodular 24-dimensional lattices.)  Our lower bounds agree
exactly with the actual numbers for $n\leq 24$. This is to be
expected for $n \leq 23$ since in those dimensions an odd or even
$n$-dimensional unimodular lattice is uniquely determined by its
root system $R$, and $\langle m(R)\, w'(R) \rangle$ is exactly
equal to the number of lattices with root system $R$ when there
are 0 or 1 such lattices. This is also true for 24-dimensional
even unimodular lattices. Our lower bound in dimension 25 is
within two percent of the actual number, and our lower bounds in
dimensions $26$ to $30$ (for which the actual numbers are not
known) are the best we are aware of.

\renewcommand{\arraystretch}{1.0} \normalsize
\begin{table}
 \caption[Lower bounds for class numbers]{Comparison of the actual
 number $\alpha_n$ of unimodular lattices
  of dimension $n$, the lower bound $\beta_n$ computed with
  equation~\ref{eqlrb}, the number $r_n$ of distinct root systems occuring
  (including those with components $\mbox{Z}$), and
   the Minkowski-Siegel mass constant $m_n.$
  Odd lattices are listed above, even lattices below.
  The last column is included only when $2 m_n$ provides
  a nontrivial lower bound for $\alpha_n$.} \label{class30}
 \begin{center}
\(
\begin{array}{crclllr}
 \dim n & \mbox{actual } \alpha_n^{\mbox{\scriptsize I \normalsize}}
 &  \geq & \mbox{bound } \beta_n^{\mbox{\scriptsize I \normalsize}}
 &  r_n^{\mbox{\scriptsize I \normalsize}} &
   \mbox{mass } m_n^{\mbox{\scriptsize I \normalsize}}
  & \beta_n^{\mbox{\scriptsize I \normalsize}}/ 2 m_n^{\mbox{\scriptsize I \normalsize}} \\
\hline

0  & 0 & \geq &0 & 0 & 0&\\

1  & 1 & \geq &1 & 1 &0.5&\\

2  & 1 & \geq & 1& 1 & 0.125&\\

3 & 1 & \geq & 1& 1 & 2.083  \times 10^{-2}&\\

4 & 1 & \geq & 1& 1 & 2.604  \times 10^{-3}&\\

5 & 1  & \geq & 1& 1 & 2.604  \times 10^{-4}&\\

6 &  1&  \geq &1& 1 & 2.170  \times 10^{-5}&\\

7 & 1& \geq &1& 1 & 1.551  \times 10^{-6}&\\

8 &  1& \geq &1& 1 & 9.688 \times 10^{-8}&\\

9 &  2& \geq &2& 2 & 6.100 \times 10^{-9}&\\

10 & 2& \geq &2& 2 & 4.485 \times 10^{-10} &\\

11 & 2& \geq &2& 2 & 4.213 \times 10^{-11}&\\

12 & 3& \geq &3& 3 & 5.267 \times 10^{-12}&\\

13 & 3& \geq &3&  3 & 9.031 \times 10^{-13}&\\

14 & 4& \geq &4& 4 & 2.186 \times 10^{-13}&\\

15 & 5& \geq &5& 5 &   7.705 \times 10^{-14}&\\

16 & 6& \geq &6&  6 & 4.093 \times 10^{-14}&\\

17 & 9& \geq &9& 9 & 3.402 \times 10^{-14}&\\

18 & 13& \geq &13& 13 & 4.583 \times 10^{-14}&\\

19 & 16& \geq &16& 16 & 1.033 \times 10^{-13}&\\

20 & 28& \geq &28& 28 & 4.002 \times 10^{-13}&\\

21 & 40& \geq &40&  40 & 2.735 \times 10^{-12}&\\

22 &  68& \geq &68& 68 & 3.377 \times 10^{-11}&\\

23 & 117 & \geq &117& 117 & 7.710 \times 10^{-10}  &\\

24 & 273 & \geq & 273& 266 & 3.330 \times 10^{-8}& \\

25  & 665 & \geq & 657 & 609 & 2.781 \times 10^{-6}&   \\

26  & ?  & \geq & 2307 & 1695 & 4.586 \times 10^{-4}& \\

27  & ?  & \geq & 14179 & 4492 & 1.524 \times 10^{-1}& \\

28  & ? &  \geq &327972 & 9213 & 1.040  \times 10^{2}& 1569.2\\

29  & ? & \geq & 37938009 & 20298 & 1.486  \times10^{5}  & 127.7
\\

30  & ?  &  \geq &20169641025 & 67848 & 4.520 \times 10^{8} &
22.3\\

31  & ? & \geq & ?  & ? &  2.980 \times 10^{12} &  \\

32  & ? & \geq & ? & ? & 4.328 \times 10^{16} &  \\

\\

 \dim n & \mbox{actual } \alpha_n^{\mbox{\scriptsize II \normalsize}}&  \geq
  & \mbox{bound } \beta_n^{\mbox{\scriptsize II \normalsize}} &
  r_n^{\mbox{\scriptsize II \normalsize}}
  &
  \mbox{mass } m_n^{\mbox{\scriptsize II \normalsize}}
  & \beta_n^{\mbox{\scriptsize II \normalsize}}/ 2 m_n^{\mbox{\scriptsize II \normalsize}} \\
\hline

 0  & 1  & \geq & 1 & 1& 1 & \\

8  & 1  & \geq & 1 &  1 & 1.435 \times 10^{-9}& \\

16  & 2 &  \geq &2  & 2 & 2.489  \times 10^{-18}& \\

24  & 24  &  \geq &24 & 24 & 7.937  \times10^{-15}  & \\

32  & ?  & \geq & 1162109024 & 13218 & 4.031  \times10^{7}  &
14.4\\

\end{array}
\)
\end{center}
\end{table}

\renewcommand{\arraystretch}{1.0} \normalsize


\section{Computing the numbers $a(R)$} \label{computeA}

For a half-integral $n \times n$ matrix $B$, define $$
c_{n,k}(B)=(-1)^{nk/2} 2^{n(k-(n-1)/2)} (\det B)^{(2k-n-1)/2} \,
b(B,k) \prod_{i=2k-n+1}^{2k}\frac{\pi^{i/2}}{\Gamma(i/2)}$$
 where $$b(B,k) = \sum_{R \in S_n(\mathrm{Q})/S_n(\mathrm{Z})}
\exp(2\pi i \,\mbox{tr}(B R))\mu(R)^{-k}$$ is the \emph{Siegel
series} (with $\mu(R)$ equal to the product of denominators of
elementary divisors of $R$). Put $\epsilon_{n,k} = 1/2$ if $n=k-1$
or $n=k>1$, and 1 otherwise.

\begin{theo}[Siegel] Let $N$ be a lattice with
$\dim(N) = n \leq 8k$, and let $$a(N) =
\frac{1}{m_{8k}^{\mbox{\scriptsize II \normalsize}}}\sum_{\Lambda
\in \Omega} \frac{r(\Lambda,N)}{|\textup{Aut}(\Lambda)|}$$ with
the sum taken over even unimodular lattices of dimension $8k$.
Then $a(N) = \epsilon_{n,8k} c_{n,4k}(B)$ where $B$ is $1/2$ times
the Gram matrix of $N$.
\end{theo}

\begin{proof}
See \cite[Theorem 6.8.1]{Kit1}, and note that the product of local
densities $\prod_p \alpha_p(M_p, N_p)$ is equal to $b(M,k)$ when
$N$ is an even unimodular lattice of dimension $2k$ \cite{Kit3}.
\end{proof}

\begin{rem} The \emph{Siegel Eisenstein series} of degree $n$ and weight $k$ is
defined to be $$E_{n,k}(Z) = \sum_{\{C,D\}} | CZ+D| ^{-k}$$ where
$\{C,D\}$ runs over all representatives of the equivalence classes
of coprime pairs of $n \times n$ matrices. For $k>n$, $c_{n,k}(C)$
is the coefficient in the Fourier expansion $$E_{n,k}(Z) =
\sum_{C}c_{n,k}(C) \exp(2\pi i \,\mbox{tr}(C Z))$$ where $C$ runs
over all positive semi-definite half-integral $n \times n$
matrices.
\end{rem}

The first explicit formula for $b(B,k)$ for arbitrary $n$ is due
to Katsurada \cite{Kat1}, and was published in 1999. Prior to
this, an explicit formula for the coefficients of $E_{n,k}(Z)$ was
known only for $n\leq 3$: The case $n=1$ is well-known (see for
example \cite[Chapter VII]{Se}), Maa\ss\ \cite{Ma} gives an
explicit formula for $n=2$ (see also Kaufhold \cite{Kau}), and
Katsurada \cite{Kat2} gives an explicit formula for $n=3$
(extending partial results by Kitaoka \cite{Kit2}).

Let $B$ be a nondegenerate symmetric half-integral $n \times n$
matrix over $\mbox{Z}_p$. It follows from \cite{Kit4} that
\small
$$b(B,s) = \left(   \zeta(s) \prod_{i=1}^{\lfloor r/2 \rfloor}
 \zeta(2s-2i)
\right) ^{-1} \prod_{p \mid D(B)} F_p(B;p^{-s}) \times \left \{
\begin{array}{ll}  L(s-r/2; \chi_B )& \mbox{if $r$ is even} \\
1 & \mbox{if $r$ is odd} \end{array} \right . $$ \normalsize
for certain polynomials $F_p(B;X)$. Here $\zeta$ is the Riemann
zeta function, and $L(s,\chi_B)$ is a Dirichlet L-series whose
values may be computed using the method in \cite{CS3}.

The explicit formula for $F_p(B;X)$ in \cite[Theorem 4.3]{Kat1} is
not itself well-suited for calculation since the outer index of
summation takes $2^n$ values, but we can use Katsurada's recursion
relations \cite[Theorems 4.1 and 4.2]{Kat1} as part of a practical
algorithm for computing $F_p(B;X)$. We will state these recursion
relations below (without proof); first we will need to introduce
some of the notation from \cite{Kat1}.

For $a = p^{r} c$ with $r \in \mathrm{Z}$ and $c \in
\mathrm{Z}_p^*$, define $\chi_p(a) = (\frac{c}{p})$ for $r$ even
and 0 for $r$ odd (where $(\frac{}{p})$ is the Legendre symbol mod
$p$), and define $$\chi_2(a) = \left \{ \begin{array}{ll} +1 &
\mbox{if $r \equiv 0 \bmod 2, \, c \equiv 1 \bmod 8$} \\ -1 &
\mbox{if $r \equiv 0 \bmod 2, \, c \equiv 5 \bmod 8$} \\ 0 &
\mbox{otherwise.}
\end{array} \right. $$
Define $\mbox{ord}_p(a)$ to be the exact power of $p$ dividing
$a$, and define $i_p(B)$ to be the least integer $t$ for which
$p^t B ^{-1}$ is half-integral. Let $(\ , \,)_p$ denote the
Hilbert Symbol over $\mathrm{Q}_p$ (see \cite{Se}), and let $h_p$
denote the Hasse invariant (see \cite{Kit1}). For odd $n$, define
$$\eta_p(B) = h_p(B) (\det B, (-1)^{(n-1)/2}\det B)_p
(-1,-1)_p^{(n^2-1)/8}.$$ For even $n$, define $$\xi_p(B) =
\chi_p((-1)^{n/2} \det B)$$ and $\xi_p
'(B)=1+\xi_p(B)-\xi_p(B)^2$. By convention, $\xi_p(B) = \xi_p'(B)
=1$ for $B$ the empty matrix. Define $D(B) = 2^{2\lfloor n/2
\rfloor}\det B$, $d_p(B) = \mbox{ord}_p(D(B))$, and $$ \delta_p
(B) =   \left \{
\begin{array}{ll} 2\lfloor (d_p(B) + 1 - \delta_{2p})/2\rfloor &
\mbox{if $n$ is even}\\ d_p(B) & \mbox{if $n$ is odd} \\
\end{array} \right. $$ where $\delta_{2p}$ is the Kronecker delta.

We shall suppress most of the subscripts $p$ in what follows. Let
$p$ be any prime and suppose $B$ and $B_2$ are nondegenerate half
integral matrices of rank $n$ and $n-1$ respectively over
$\mathrm{Z}_p$. Put $\delta = \delta(B)$ and $\tilde{\delta} =
\delta(B_2)$. If $n$ is even, put  $\xi = \xi(B)$, $\xi' =
\xi'(B)$, and $\tilde{\eta} = \eta(B_2)$; if $n$ is odd, put
$\tilde{\xi} = \xi(B_2)$, $\tilde{\xi}' = \xi'(B_2)$, and $\eta =
\eta(B)$. (By convention $\tilde{\eta} = \tilde{\eta}'=1$ if $B_2
= \emptyset$, and $\tilde{\delta}=0$ if $n=1$.) Then define
rational functions $C(B,B_2;X)^{(1)}$ and  $C(B,B_2;X)^{(0)}$ in
$X$ by
   $$C(B,B_2;X)^{(1)} = \left \{  \begin{array}{ll}
   \frac{1-p^{n/2} \xi X}{1 - p^{n+1}X^2} & \mbox{if $n$ is even} \\
    \frac{1}{1 - p^{(n+1)/2} \tilde{\xi} X} & \mbox{if $n$ is odd} \end{array}
\right. $$ and
 $$C(B,B_2;X)^{(0)} = \left \{
\begin{array}{ll}
   \frac{(-1)^{\xi + 1} \xi' \tilde{\eta}(1-p^{n/2+1}X \xi)
   (p^{n/2}X)^{\delta - \tilde{\delta}+\xi^2} p^{\delta /2})}
   {1 - p^{n+1}X^2} & \mbox{if $n$ is even} \\
    \frac{(-1)^{\tilde{\xi}} \tilde{\xi'} \eta
    (p^{(n-1)/2}X)^{\delta - \tilde{\delta} + 2 - \tilde{\xi}^2}
    p^{(2 \delta - \tilde{\delta} + 2)/2}}
    {1 - p^{(n+1)/2} \tilde{\xi} X} & \mbox{if $n$ is odd.} \end{array}
\right. $$

\begin{theo}[Katsurada \cite{Kat1}] \label{katp} Let $B_1 = (b_1)$ and $B_2$ be
nondegenerate half integral matrices of degree $1$ and $n-1$,
respectively, over $\mathrm{Z}_p$, and put $B=B_1 \bot B_2$.
Assume that $ord(b_1) \geq i(B_2) - 1 + 2\delta_{2p}$. Then we
have $$F_p(B;X) = C(B,B_2;X)^{(1)} F(B_2;pX) +
C(B,B_2;X)^{(0)}F_p(B_2;X).$$
\end{theo}


Let $B_2$ be a nondegenerate half-integral matrix of degree $n-2$
over $\mathrm{Z}_2$, let $H = {0 \ \ 1/2 \choose 1/2 \ \ 0} $ and
let $Y = {1 \ \ 1/2 \choose 1/2 \ \ 1}$.
Let $B_1 = 2^m K$ with $K = H$ or $Y$, or $B_1=  2^m u_1 \bot 2^m
u_2$ with $ u_1, u_2 \in \mathrm{Z}_2^{*}$, and put $B$ = $B_1
\bot B_2$, $\delta = \delta(B)$, $\tilde{\delta} = \delta(2^m \bot
B_2)$, $\hat{\delta} = \delta(B_2)$, and $$ \sigma =   \left \{
\begin{array}{ll} (2 \tilde{\delta} - \delta - \hat{\delta}+2)/2 &
\mbox{if $n$ is even, $B_1 = 2^m u_1 \bot 2^m u_2$ and $d(B)$ is
odd}\\ & \mbox{or if $n$ is even, $B_1 = 2^m K$, and $\xi(B_2)=0$}
\\ 2 & \mbox{if $n$ is odd, $B_1 = 2^m K$, and $d(2^m \bot B_2)$ is
even} \\ 0 & \mbox{otherwise.} \\ \end{array} \right. $$ If $n$ is
even, put $\xi = \xi(B)$, $\xi' = \xi'(B)$, $\hat{\xi} =
\xi(B_2)$, $\hat{\xi}' = \xi'(B_2)$, and  $$ \tilde{\eta} = \left
\{ \begin{array}{ll} \eta (2^m u_2 \bot B_2) & \mbox{if $B_1 = 2^m
u_1 \bot 2^m u_2$,} \\ & \mbox{and $d(B_2)$ is even}
\\ (-1)^{((n - 1)^2 - 1)/8} h(B_2)(2^m, (-1)^{(n-1)/2} \det B_2)_2
& \mbox{if $B_1 = 2^m K$,} \\ & \mbox{and $\xi (B_2) \neq 0$}
\\ 1 & \mbox{otherwise.} \\ \end{array} \right. $$
If $n$ is odd, put $\eta$ = $\eta(B)$, $\hat{\eta} = \eta(B_2)$,
$\tilde{\xi}' = 1$, and $\tilde{\xi} = 1 $ if $B_1 = 2^m K$ and
$d(2^m \bot B_2)$ is even, and 0 otherwise. Define four rational
functions in $X$ by
 $$C(B,B_2;X)^{(11)} =
\left \{  \begin{array}{ll}
   \frac{1-2^{n/2} \xi X}{1 - 2^{n+1}X^2} & \mbox{if $n$ is even} \\
    \frac{1}{1 - 2^{(n+1)/2} \tilde{\xi} X} & \mbox{if $n$ is odd} \end{array}
\right. $$ $$C(B,B_2;X)^{(10)} = \left \{
\begin{array}{ll}
   \frac{(-1)^{\xi + 1} \xi' \tilde{\eta}(1-2^{n/2+1}X \xi)
   (2^{n/2}X)^{\delta - \tilde{\delta}+\xi^2 + \sigma} 2^{\delta /2})}
   {1 - 2^{n+1}X^2} & \mbox{if $n$ is even} \\
    \frac{(-1)^{\tilde{\xi}} \eta
    (2^{(n-1)/2}X)^{\delta - \tilde{\delta} + 2 - \tilde{\xi}^2 + \sigma}
    2^{(2 \delta - \tilde{\delta} + 2 + \sigma)/2}}
    {1 - 2^{(n+1)/2} \tilde{\xi} X} & \mbox{if $n$ is odd} \end{array}
\right. $$
 $$C(B,B_2;X)^{(21)} = \left \{  \begin{array}{ll}
 \frac{1}{1 - 2^{n/2} \hat{\xi} X} & \mbox{if $n$ is even} \\
   \frac{1-2^{(n-1)/2} \tilde{\xi} X}{1 - 2^n X^2} & \mbox{if $n$ is odd} \end{array}
\right. $$ $$C(B,B_2;X)^{(20)} = \left \{
\begin{array}{ll}
  \frac{(-1)^{\hat{\xi}} \hat{\xi'} \tilde{\eta}
    (2^{(n-2)/2}X)^{\tilde{\delta} - \hat{\delta} + 2 - \hat{\xi}^2 - \sigma}
    2^{(2 \tilde{\delta} - \hat{\delta} + 2 - 2\sigma)/2}}
    {1 - 2^{n/2} \hat{\xi} X} & \mbox{if $n$ is even} \\
   \frac{(-1)^{\tilde{\xi} + 1} \hat{\eta}(1-2^{(n+1)/2}X \tilde{\xi})
   (2^{(n-1)/2}X)^{\tilde{\delta} - \hat{\delta}+\tilde{\xi}^2 - \sigma}
    2^{(\tilde{\delta} - \sigma) /2})}
   {1 - 2^n X^2} & \mbox{if $n$ is odd.}
   \end{array}
\right. $$

\begin{theo}[Katsurada \cite{Kat1}] \label{kat2}
Let $B_1 = 2^m u_1 \bot 2^m u_2$ with $u_1, u_2 \in
\mathrm{Z}_2^{*}$ or $B = 2^m K$ with $K = H$ or $Y$. Let $B_2$ be
a half-integral matrix of degree $n-2$ over $\mathrm{Z}_2$ which
is also in $GL_{n-2}(\mathrm{Q}_2)$, and put $B=B_1 \bot B_2$.
Assume that $m \geq i(B_2) + 1$. Then we have
\begin{eqnarray*}
F_2(B;X) & = & C(B,B_2;X)^{(11)}C(B,B_2;2X)^{(21)} F_2(B_2;4X) \\
 & + &  C(B,B_2;X)^{(11)} C(B,B_2;2X)^{(20)} F_2(B_2;2X) \\
 & + &  C(B,B_2;X)^{(10)} C(B,B_2;X)^{(21)} F_2(B_2;2X) \\
 & + &  C(B,B_2;X)^{(10)} C(B,B_2;X)^{(20)} F_2(B_2;X).
\end{eqnarray*}
\end{theo}

For $p \neq 2$, any nondegenerate symmetric half-integral $n
\times n$ matrix $B$ can be diagonalized over the $p$-adic
integers $\mbox{Z}_p$, $ B \cong p^{e_1}u_1 \bot \ldots \bot
p^{e_n}u_n $ with $e_1 \geq \cdots \geq e_n \geq 0$ and with $u_i
\in \{ 1,\epsilon \}$ for all $i$, where $\epsilon$ is any
quadratic nonresidue (see \cite{Wa}). Then $F_p(B;X)$ can be
computed by repeated applications of Theorem~\ref{katp} above.

Note that in the course of recursively computing $F_p(B;X)$ for
$X=p^{-k}$ in this manner, one occasionally encounters zeros of
the denominators of the functions $C^{(1)}$ and $C^{(0)}$. (This
is only a problem when $X$ is a negative power of $p$.) Rather
than attempting to simplify the expressions symbolically, we
instead computed $F_p(B;1), F_p(B;p), \ldots ,F_p(B;p^n)$, and
then used Lagrangian interpolation to compute $F_p(B;p^{-k})$. The
reason for choosing $X=1, p, \ldots , p^n$ rather than $n+1$ other
numbers is that for these numbers the recursive subproblems
overlap: we need only evaluate $F_p(p^{e_i}u_i \bot \ldots \bot
p^{e_n}u_n; X)$ for $X = 1,p, \ldots, p^{2n-i}$, starting with
$i=n$, and working down to $i=1$ which gives $F_p(B;1), \ldots
,F_p(B;p^n)$.

For $p=2$, any nondegenerate symmetric half-integral matrix $B$ is
equivalent over $\mbox{Z}_2$ to a matrix of the form $2^{e_1}(U_1
\bot V_1) \bot \ldots \bot 2^{e_m}(U_m \bot V_m) $ with $e_1 >
\cdots > e_m \geq 0$, $U_i = \emptyset$ or $u_1$ or $u_1 \bot u_2$
for $u_1,u_2 \in \{\pm 1, \pm3 \} $, and $V_i = \emptyset$ or $H
\bot \cdots \bot H$ or $H \bot \cdots \bot H \bot Y $ (see
\cite{Wa}). Then $F_2(B;X)$ can be computed by repeated
applications of Theorems~\ref{katp} and \ref{kat2} above. (We
again used an interpolation scheme, similar to the one described
for $p \neq 2$.)

 \section{Computing the numbers $r(R_i, R_j)$} \label{computeR}

The method described in this section is essentially the one used
for some of the computations in \cite{BFW} (although the algorithm
itself is not described in that paper). In this section we write
$r(R,R')$ as $\mbox{emb}(R', R)$, since a representation of $R'$
by $R$ is the same as a linear map from $R'$ into $R$ which
preserves inner products. For any irreducible root systems $S$ and
$T$ it is routine to compute the number of embeddings
$\mbox{emb}(S, T)$ of $S$ into $T$, and to determine the root
system of the orthogonal complement of $S$ in $T$ for each of
these embeddings. There are at most two orbits of embeddings of
$S$ into $T$, so we write $\mbox{emb}(S, T) = \mbox{emb}_1(S, T) +
\mbox{emb}_2(S, T)$, where there are $\mbox{emb}_1(S,T)$
embeddings of S into T for which the orthogonal complement is
$\mbox{comp}_1(S,T)$ and there are $\mbox{emb}_2(S,T)$ embeddings
of $S$ into $T$ for which the orthogonal complement is
$\mbox{comp}_2(S,T)$.

\begin{exa}
We demonstrate how to compute the two orbits of embeddings of
$A_3$ into $D_n$. Recall from \cite{CS1} that $A_n = \{(x_1,
\ldots, x_{n+1}) \in \mbox{Z}^{n+1} : x_1 + \cdots + x_{n+1} =
0\}$ and $D_n = \{(x_1, \ldots, x_n) \in \mbox{Z}^{n} : x_1 +
\cdots + x_n \in 2\mbox{Z}\}$. $A_3$ is generated by three roots
$v_1=(1,-1,0,0)$, $v_2 = (0,1,-1,0)$, and $v_3 = (0,0,1,-1)$ with
$(v_1,v_1) =(v_2,v_2) =(v_3,v_3) = 2$, $(v_1,v_2) = (v_2,v_3)=-1$
and $(v_1,v_3)=0$. We can map $v_1$ to any of the $4 C(n,2)$ roots
$r_1$ in $D_n$ (all permutations of $(\pm 1,\pm 1,0, \ldots ,
0)$), where $C(n,k)$ denotes the binomial coefficient.

We can then map $v_2$ to any of the $4(n-2)$ roots $r_2$ of $D_n$
which have inner product $-1$ with $r_1$.  Let us say $r_1$ is
supported in coordinates $i$ and $j$ and $r_2$ is supported in
coordinates $j$ and $k$. (Clearly $k \neq i$.) Then there are two
cases:

(a) We can map $v_3$ to the root $r_3$ supported in coordinates
$i$ and $j$ which has inner product $0$ with $r_1$ and inner
product $-1$ with $r_2$. In this case the roots of $D_n$
orthogonal to $r_1$, $r_2$, and $r_3$ form the system $D_{n-3}$;

(b) We can map $v_3$ to $2(n-3)$ roots $r_3$ supported in
coordinates $k$ and $l \,\,(l \neq i,j)$ which have inner product
0 with $r_1$ and inner product $-1$ with $r_2$. In this case the
roots of $D_n$ orthogonal to $r_1$, $r_2$, and $r_3$ form the
system $D_{n-4}$.

Hence there are $4 C(n,2) \cdot 4(n-2) \cdot 1 = C(n,3) \cdot
|\mbox{Aut}(A_3)|$ ways to embed $A_3$ into $D_n$ with complement
$D_{n-3}$, and there are  $4 C(n,2) \cdot 4(n-2) \cdot 2(n-3) =
2^{3} C(n,4)\cdot  |\mbox{Aut}(A_3)|$ ways to embed $A_3$ into
$D_n$ with complement $D_{n-4}$.
\end{exa}

The computations for other irreducible root systems are similar.
All nonzero values of $\mbox{emb}_1$, $\mbox{emb}_2$,
$\mbox{comp}_1$, and $\mbox{comp}_2$ are given in
Table~\ref{embeddings}.

\renewcommand{\arraystretch}{1.0} \normalsize

\begin{table}
 \caption[Embeddings of irreducible root systems into one another]
 {Embeddings of irreducible root systems into one another other. (In the fourth and last
columns, $A_0, D_0, D_1, D_2 \mbox{ and } D_3$ should be
interpreted as $\emptyset, \emptyset, \emptyset, A_1 A_1$ and
$A_3$, respectively.)} \label{embeddings}
\begin{center}

\small

\(
\begin{array}{cccccc}
S  & T & \frac{\mathrm{emb}_1(S,T)}{|\mathrm{Aut}(S)|} &
\mbox{comp}_1(S,T)& \frac{\mathrm{emb}_2(S,T)}{|\mathrm{Aut}(S)|}
& \mbox{comp}_2(S,T)
\\ \hline

 \emptyset & T & 1 & T & - & - \\

A_i & A_j \,(j \geq i) & C(j+1, i+1) & A_{j-i-1} &- &- \\

A_1 & D_j & C(j,2)\cdot 2 & A_1 D_{j-2} & -& - \\

A_3 & D_j & C(j,4)\cdot 2^3 & D_{j-4} & C(j,3) & D_{j-3} \\

A_i\, (i \neq 1,3)& D_j\, (j > i) & C(j,i+1)\cdot 2^i & D_{j-i-1}
&- &-
\\

D_i & D_j\, ( j \geq i) & C(j,i) & D_{j-i}  & -&-  \\

A_1 & E_6 &  2^2 3^2&  A_5  & -&- \\

A_2 & E_6  &  2^3 3 \cdot 5 &  A_2  A_2  &- &- \\

A_3 & E_6 &  2 \cdot 3^3 5 & A_1  A_1  &- &-\\

A_4 & E_6 &  2^3 3^3 & A_1  & -&- \\

A_5 & E_6 &  2^2 3^2 &  A_1   & -&-\\

D_4 & E_6 & 3^2 5 &  \emptyset  & -&-\\

D_5 & E_6  & 3^3  &   \emptyset   & -&- \\

E_6 & E_6  & 1  &  \emptyset   & -&- \\

A_1 & E_7 &  3^2 7&  D_6   & -&-\\

A_2 & E_7 & 2^4  3 \cdot 7 &  A_5 &- &- \\

A_3 & E_7 & 2^2 3^2  5 \cdot 7 &  A_3 A_1 &- &-  \\

A_4 & E_7  & 2^5 3^2 7 &  A_2 & -&- \\

A_5 & E_7  & 2^4  3 \cdot 7 & A_2   & 2^4 3^2 7 & A_1 \\

A_6 & E_7  &  2^5 3^2 &  \emptyset &- &-  \\

A_7 & E_7  &  2^2 3^2 &   \emptyset & -&-  \\

D_4 & E_7 &  3^2 5 \cdot 7& A_1 A_1 A_1 &- &-  \\

D_5 & E_7 &  2 \cdot 3^3 7& A_1 &- &- \\

D_6 & E_7 &  3^2 7 & A_1  & -&- \\

E_6 & E_7 &  2^2 7 &  \emptyset &- &-  \\

E_7 & E_7 &  1 &  \emptyset  & -&- \\

A_1 & E_8  & 2^3 3 \cdot 5 & E_7 &- &-  \\

A_2 & E_8  & 2^5 5 \cdot 7 &  E_6 & -& -\\

A_3 & E_8  & 2^3 3^3 5 \cdot 7 & D_5  & -&- \\

A_4 & E_8  & 2^7 3^3 7 & A_4 & -&- \\

A_5 & E_8  & 2^7 3^2 5 \cdot 7  & A_2 A_1 &- &- \\

A_6 & E_8  & 2^8 3^3 5 &  A_1 & -&- \\

A_7 & E_8  & 2^5 3^3 5 &  A_1  & 2^6 3^3 5 & \emptyset\\

A_8 & E_8  & 2^6 3 \cdot 5 & \emptyset  & -&-\\

D_4 & E_8  & 2 \cdot 3^2 5^2 7 & D_4 & -&- \\

D_5 & E_8  & 2^3 3^3 5 \cdot 7 & A_3  & -&-\\

D_6 & E_8  & 2^2 3^3 5 \cdot 7 &  A_1 A_1 & -&- \\

D_7 & E_8  & 2^3 3^3 5 &   \emptyset & -&-  \\

D_8 & E_8  & 3^3 5 & \emptyset  & -&- \\

E_6 & E_8  & 2^5 5 \cdot 7 & A_2 & -&-  \\

E_7 & E_8  & 2^3 3 \cdot 5 & A_1 & -&-\\

E_8 & E_8  & 1 &   \emptyset  & -&- \\

\end{array}
\)
\end{center}
\normalsize \end{table}

\renewcommand{\arraystretch}{1.0} \normalsize

 If $R$ and $R'$ are root
lattices then we can write $R = S_1 \oplus \cdots \oplus S_k$ and
$R' = T_1 \oplus \cdots \oplus T_m$ where each $S_i$ and $T_i$ is
an irreducible root lattice, and we can compute $r(R',R)$
recursively via the formula
\begin{eqnarray*}
r(R', R) & = & \mbox{emb}(R,R') \\ &  = & \mbox{emb}
(\bigoplus_{i\leq k} S_i , \bigoplus_{i \leq m} T_i) \\ & = &
\sum_{j=1}^m \mbox{emb}_1(S_k, T_j) \, \mbox{emb}(\bigoplus_{i
\leq k-1}S_i, \mbox{comp}_1(S_k, T_j) \oplus \bigoplus_{i\leq m, i
\neq j}T_i)\\ & + & \sum_{j=1}^m \mbox{emb}_2(S_k, T_j) \,
\mbox{emb}(\bigoplus_{i \leq k-1}S_i, \mbox{comp}_2(S_k, T_j)
\oplus \bigoplus_{i\leq m, i \neq j} T_i ).
\end{eqnarray*}

A direct implementation of this algorithm does a lot of redundant
computation on certain inputs, some of which we can circumvent
with dynamic programming or \emph{memoization} (see \cite[Chapter
16]{CLR}). Our implementation also does several things to reduce
the amount of computation when there are direct summands in $R$ or
$R'$ with multiplicity greater than one. But computing $r(R',R)$
is an NP-Hard problem, since it is an NP-Complete problem to
determine if $r(R',R)>0$, as we show below.

\begin{prop}
The problem of determining whether a root system $R$ embeds into a
root system $R'$ is NP-Complete.
\end{prop}


\begin{proof}
We shall reduce \textsc{3-Partition} (see \cite{GJ}) to this
problem. Let $S = \{s_1, \ldots, s_{3k}\}$ be an instance of
\textsc{3-Partition}, with the $s_i$ positive integers summing to
$kt$. Then $S$ can be partitioned into $k$ sets each consisting of
3 elements with sum $t$ if and only if the root system $D_{4s_1}
\oplus \cdots \oplus D_{4s_{3k}}$ embeds into the root system
$\bigoplus_{i=1}^{k}D_{4t}$.
\end{proof}

\begin{rem} Since \textsc{3-Partition} is strongly NP-Complete,
the problem of whe\-ther one root system embeds into another remains
NP-Complete if the Gram matrices of root lattices are used as
input (rather that the list of components $A_i$, $D_i$, and $E_i$
as above).
\end{rem}

Since we need to compute $r(R, R')$ for all pairs of root systems
rather than just one pair, the amortized computational cost would
be reduced considerably by using dynamic programming. But the
dynamic programming table becomes unmanageably large in dimension
32, so we instead use a hash table, which is purged periodically,
for memoization.

\section{Eliminating root systems \emph{a priori}} \label{hacks}

There are 405844 root systems of rank $n \leq 32$ with no vectors
of norm 1, corresponding to all direct sums of $A_i, D_i\,(i \geq
4), E_6, E_7 \mbox{ and } E_8$, where the order of the summands
does not matter and the sum of the subscripts is at most 32. Since
computing the number of embeddings $r(R_i, R_j)$ can be
time-consuming, we used the following congruences, due to
Borcherds \cite{Bo1}, to eliminate some root systems from
consideration:

Let $\mbox{roots}(R)$ denote the number of roots of $R$. If $R$ is
the root system of a 32-dimensional even unimodular lattice then
\begin{center}
\begin{quote}

if $R$ contains $E_8$, then $\mbox{roots}(R) \equiv 0 \pmod
  {24}$ \\ \noindent
if $R$ contains $E_7$, then $\mbox{roots}(R) \equiv 0 \pmod
  {12}$\\ \noindent
 if $R$ contains $E_6$, then $\mbox{roots}(R)
\equiv 0 \pmod
  {6}$\\ \noindent
  if $R$ contains $D_6$, then $\mbox{roots}(R)
\equiv 0 \pmod
  {4}$\\ \noindent
  if $R$ contains $D_7$, then $\mbox{roots}(R)
\equiv 0 \pmod
  {8}$\\ \noindent
  if $R$ contains $D_8$, then $\mbox{roots}(R)
\equiv 0 \pmod
  {8}$\\ \noindent
 if $R$ contains $D_n, n>8$, then $\mbox{roots}(R)
\equiv 0 \pmod
  {16}$\\
\end{quote}
\end{center}

Also note that if $a(R)=0$ then $m(R)=0$. If $\mbox{rank}(R)=32$
and $\det(R)$ is not a perfect square then $a(R)$ must be 0, so we
eliminated those root systems as well. This left 135443 root
systems. We ordered them so that $\dim(R_i)\leq \dim(R_j)$ if
$i<j$, and so that $\det(R_i) \geq \det(R_j)$ if $i<j$ and
$\dim(R_i) = \dim(R_j)$.

As the matrix $U$ with $U_{i,j} = r(R_i,R_j)$ would still contain
around 10 billion elements, we did not explicitly construct and
invert it. Rather, we computed each element in the matrix when it
was required for solving $\frac{1}{m}U v = w$ by
back-substitution, with  $$m(R_i) = \frac{1}{r(R_i,R_i)}\left\{m
\cdot a(R_i) - \sum_{j>i}r(R_j, R_i)\, m(R_j)\right\} .$$

If $m(R_j)$ has already been computed to be 0, then the values
$r(R_j, R_i)$ need not be computed since they make no contribution
to this sum.

\section*{Acknowledgements}
The author would like to thank Richard Borcherds for many helpful
suggestions, Richard Fateman for advice on Lisp, and Xerox PARC
for the use of its computers. This work was partially supported by
grants from the NSF and the Royal Society.

\newpage

\bibliographystyle{plain}

\end{document}